\documentclass{IEEEtran}
                                                                                                                                                                                                                                 
\pdfoutput=1
\usepackage{mathtools}
\usepackage{graphicx}
\usepackage{algorithm}
\usepackage{algpseudocode}
\usepackage{amsmath}
\usepackage{amsfonts}
\usepackage{verbatim}
\usepackage{framed}
\usepackage{todonotes}

\usepackage{float}
\floatstyle{boxed} 
\restylefloat{algorithm}

\usepackage{eqparbox}

\newcommand{\half}{\frac{1}{2}}

\newcommand{\eqb}[1]{\begin{equation}\label{#1}}
\newcommand{\eqe}{\end{equation}}
\newcommand{\splitb}{\begin{align}\begin{split}}
\newcommand{\splite}{\end{split} \end{align}}

\newcommand{\st}{\hbox{ subject to }}
\newcommand{\matb}{\left( \begin{matrix*}[r] }
\newcommand{\mate}{\end{matrix*}\right)}

\newcommand{\kp}{_{k+1}}
\newcommand{\opt}{^\star}
\newcommand{\reals}{\mathbb{R}}

\newcommand{\mot}{\frac{\mu}{2}}
\newcommand{\oott}{\frac{1}{2\tau}}
\newcommand{\oots}{\frac{1}{2\sigma}}
\newcommand{\tinv}{\frac{1}{\tau}}

\newcommand{\tv}{|\nabla x|}
\newcommand{\tvmax}{y\cdot \nabla x}
\newcommand{\A}{\textbf{A}}
\newcommand{\B}{\textbf{B}}
\newcommand{\C}{\textbf{C}}
\newcommand{\Cone}{\textbf{C1}}
\newcommand{\Ctwo}{\textbf{C2}}

\DeclareMathOperator*{\argmin}{arg\,min}

\newtheorem{theorem}{Theorem}
\newtheorem{lemma}{Lemma}
\newtheorem{remark}{Remark}

\begin{document}

\sloppy

\bibliographystyle{ieeetr} 
\title{Adaptive Primal-Dual Hybrid Gradient Methods for Saddle-Point Problems}
\author{Tom Goldstein, Min Li, Xiaoming Yuan, Ernie Esser, Richard Baraniuk}
\date{\today}
\maketitle

\begin{abstract}
The Primal-Dual hybrid gradient (PDHG) method is a powerful optimization scheme that breaks complex problems into simple sub-steps. Unfortunately, PDHG methods require the user to choose stepsize parameters, and the speed of convergence is highly sensitive to this choice.  We introduce new adaptive PDHG schemes that automatically tune the stepsize parameters for fast convergence without user inputs.  We prove rigorous convergence results for our methods, and identify the conditions required for convergence.  We also develop practical implementations of adaptive schemes that formally satisfy the convergence requirements.  Numerical experiments show that adaptive PDHG methods have advantages over non-adaptive implementations in terms of both efficiency and simplicity for the user.
\end{abstract}



\section{Introduction}
This manuscript considers saddle-point problems of form
\eqb{saddle}
\min_{x\in X} \max_{y\in Y} f(x)+ y^T Ax - g(y)
\eqe
where $f$ and $g$ are convex functions, $A\in \reals^{M\times N}$ is a matrix, and $X\subset \reals^N$ and $Y\subset \reals^M$ are  convex sets.

 The formulation (\ref{saddle}) is particularly appropriate for solving inverse problems involving the $\ell_1$ norm. Such problems take the form
    \eqb{l1_gen}\min_x | S x| + \frac{\mu}{2}\|Bx-f \|^2\eqe 
     where $|\cdot|$ denotes the $\ell_1$ norm, $B$ and $S$  are linear operators, and $\|\cdot\|$ is the $\ell_2$ norm.  The formulation \eqref{l1_gen} is useful because it naturally enforces sparsity of $Sx.$ Many different problems can be addressed by choosing different values for the sparsity transform $S.$ 
          
          In the context of image processing, $S$ is frequently the gradient operator.
  In this case the $\ell_1$ term becomes  $|\nabla u|,$ the total variation semi-norm \cite{ROF92}.    Problems of this form arise whenever an image is recovered from incomplete or noise-contaminated data.  In this case, $B$ is often a Gaussian blur matrix or a sub-sampled fast transform, such as a Fourier or Hadamard transform.

As we will see below, the problem \eqref{l1_gen} can be put in the ``saddle-point'' form \eqref{saddle} and can then be addressed using the techniques described below.  Many other common minimization problems can also be put in this common form, including image segmentation, TVL1 minimization, and general linear programing.  

 In many problems of practical interest, $f$
 and $g$ do not share common properties, making it difficult to
 derive numerical schemes for (\ref{saddle}) that address both
 terms simultaneously.   However, it frequently occurs in practice that
 efficient algorithms exist for minimizing $f$ and $g$  independently.  In this
 case, the Primal-Dual Hybrid Gradient (PDHG)  \cite{EZC09,CP10} method is quite useful.   This method removes the coupling between $f$ and $g,$ enabling each term to be addressed separately.   Because it decouples $f$ and $g$, the steps of PDHG can often be written explicitly, as opposed to other splitting methods that require expensive minimization sub-problems.

One of the primary difficulties with PDHG is that it relies on step-size parameters that must be carefully chosen by the user.  The speed of the method depends heavily on the choice of these parameters, and there is often no intuitive way to choose them.  Furthermore, the stepsize restriction for these methods depends on the spectral properties of $A$, which may be unknown for general problems.

 We will present practical adaptive schemes that optimize the PDHG parameters automatically as the problem is solved.  Our new methods are not only much easier to use in practice, but also result in much faster convergence than constant-stepsize schemes.   After introducing the adaptive methods, we prove new theoretical results that guarantee convergence of PDHG under very general circumstances, including adaptive stepsizes.

\subsection{Notation}
Given two vectors $u,v\in \reals^N,$ we will denote their discrete inner product by $u^T v  = \sum_i u_i v_i.$   When $u,v\in \reals^\Omega$ are not column vectors, we will also use the ``dot'' notation for the inner product, $u \cdot v = \sum_{i\in \Omega} u_iv_i.$   \vspace{2pt}

The formulation \eqref{saddle} can be generalized to handle complex-valued vectors.  In this case, we will consider the real part of the inner product.  We use the notion $\Re\{\cdot\}$ to denote the real part of a vector or scalar.

  We will make use of a variety of norms, including the  $\ell_2$-norm, $ \|u\| = \sqrt{\sum_i u_i^2}$, and the $\ell_1$-norm, $ |u| = \sum_i |u_i|$.\

When $M$ is a  symmetric positive definite matrix, we define the $M$-norm by 
$$\|u\|_M^2 =  u^T Mu.$$
If $M$ is indefinite, this does not define a proper norm.  In this case, we will still  write $\|u\|_M^2$ to denote the quantity $u^T Mu,$ even though this is an abuse of notation.

  For a matrix M, we can also define the ``operator norm''
  $$\|M\|_{op} = \max_u \frac{\| Mu\|}{\|u\|}.$$
  If $M$ is symmetric, then the operator norm is the spectral radius of $M,$ which we denote by $\rho(M).$

We use the notation $\partial f$ to denote the sub-differential (i.e., generalized derivative) of a function $f$.   

Finally, we will use $\chi_C$ to denote the characteristic function of a convex set $C,$ which is defined as follows:  
 $$
 \chi_C(x) = 
 \begin{cases}
 0, &\text{ if } x\in C
 \\ \infty, &\text{  otherwise}.
 \end{cases}
 $$

\section{The Primal-Dual Hybrid Gradient Method}
The PDHG scheme has its roots in the well-known Arrow-Hurwicz scheme, which was originally proposed in the field of economics and later refined for solving saddle point problems by Popov \cite{Popov80}.  While the simple structure of the Arrow-Hurwicz method made it appealing, tight stepsize restrictions and poor performance make this method impractical for many problems.

Research in this direction was reinvigorated by the introduction of PDHG, which converges rapidly for a wide range of stepsizes.   PDHG originally appeared in a technical report by Zhu and Chan \cite{ZC08}.  It was later analyzed for convergence in \cite{EZC09,CP10}, and studied in the context of image segmentation in \cite{PCBC09:seg}.  An extensive technical study of the method and its variants is given by He and Yuan \cite{HY12}.  Several extensions of PDHG, including simplified iterations for the case that $f$ or $g$ is differentiable, are presented by Condat \cite{Condat13}.

PDHG is listed in Algorithm \ref{alg:pdhg}.   The algorithm is inspired by the forward-backward algorithm, as the primal and dual parameters are updated using a combination of forward and backward steps. In steps (2-3), the method updates $x$ to decreases the energy \eqref{saddle} by first taking a gradient descent step with respect to the inner product term in \eqref{saddle}, and then taking a ``backward'' or proximal step involving $f$. In steps (5-6), the energy \eqref{saddle} is increased by first marching up the gradient of the inner product term with respect to $y$,  and then a backward step is taken with respect to $g$. 

\begin{algorithm}
    \caption{Basic PDHG}
    \label{alg:pdhg}
    \begin{algorithmic}[1]
    \Require $x_0\in \reals^N$, $y_0\in \reals^M$, $\sigma_k,\tau_k>0$
    \While{\emph{Not Converged}}
	\State $\hat x\kp = x_k - \tau_k A^Ty_k$  
    \State $x\kp = \argmin_{x\in X} f(x)+\frac{1}{2\tau_k}\|x-\hat x\kp \|^2  $ \label{proxF}
    \State $\bar x\kp = x\kp+(x\kp-x_k)$ \label{pdhg:predict}
     \State $\hat y\kp = y_k + \sigma_k A\bar x\kp$
     \State $y\kp = \argmin_{y\in Y}  g(y) +\frac{1}{2\sigma_k}\|y-\hat y\kp \|^2 $  \label{proxG}   
    \EndWhile
       \end{algorithmic}
  \end{algorithm}

Note that steps \ref{proxF} and \ref{proxG} of Algorithm \ref{alg:pdhg} require minimizations.  These minimization steps can be written in a more compact form using the \emph{proximal} operators of $f$ and $g$:
  \begin{align} \begin{split}
  J_{\tau F}(\hat x) &= \argmin_{x \in X} f(x)+\frac{1}{2\tau}\|x-\hat x \|^2 = (I+\tau F)^{-1} \hat x\\ 
  J_{\sigma G}(\hat y) &= \argmin_{y \in Y} g(y)+\frac{1}{2\sigma}\|y-\hat y \|^2 = (I+\sigma G)^{-1} \hat y.
 \label{resolvents}
 \end{split} \end{align}

Algorithm \ref{alg:pdhg} has been analyzed in the case of constant stepsizes, $\tau_k=\tau$ and $\sigma_k=\sigma.$  In particular, it is known to converge as long as  $\sigma\tau <\frac{1}{\rho(A^TA)}$ \cite{EZC09, CP10, HY12}.  However, PDHG typically does not converge when non-constant stepsizes are used, even in the case that $\sigma_k\tau_k <\frac{1}{\rho(A^TA)}.$

In this article, we identify the specific stepsize conditions that guarantee convergence and propose practical adaptive methods that enforce these conditions.

\begin{remark}
Step 4 of Algorithm \ref{alg:pdhg} is a prediction step of the form $\bar x\kp = x\kp+\theta(x\kp-x_k),$ where $\theta=1.$  Note that PDHG methods have been analyzed for $\theta \in [-1,1]$ (see \cite{HY12}). Some authors have even suggested non-constant values of $\theta$ as a means of accelerating the method  \cite{CP10}.  However, the case $\theta=1$ has been used almost exclusively in applications.   
\end{remark}

\subsection{Optimality Conditions and Residuals}
Problem \eqref{saddle} can be written in its unconstrained form using the characteristic functions for the sets $X$ and $Y$.  We have  
\eqb{saddleInd}
\min_{x\in \reals^N} \max_{y\in \reals^M} f(x)+\chi_X(x) +y^T Ax - g(y)-\chi_Y(y) .
\eqe
 	Let $F = \partial \{f+\chi_X(x) \} ,$ and $G=\partial \{g+\chi_Y(y)\}$\footnote{Note that $\partial \{f+\chi_X(x) \}=\partial f+\partial\chi_X(x) $ precisely when the resolvent minimizations \eqref{resolvents} are feasible (see Rockafellar \cite{Rockafellar96} Theorem 23.8), and in this case the optimality conditions \eqref{saddleInd} admit a solution.  We will assume for the remainder of this article that the problem \eqref{saddleInd} is feasible. }.  From the equivalence between  \eqref{saddle} and \eqref{saddleInd}, we see that the optimality conditions for \eqref{saddle} are
\begin{align} \label{opt1} 
0 &\in F(x\opt)+A^Ty\opt  \\
0 &\in G(y\opt) - A(x\opt) .  \label{opt2} 
 \end{align}
	The optimality conditions (\ref{opt1},\ref{opt2}) state that the derivative with respect to both the primal and dual variables must be zero.
This motivates us to define the primal and dual residuals

\begin{align}\begin{split} \label{resid} 
P(x,y) &= F(x)+A^Ty \\
D(x,y) &= G(y) - Ax.  
\end{split} \end{align}
We can measure convergence of the algorithm by tracking the size of these residuals.  Note that the residuals are in general multi-valued (because they depend on the subdifferentials of $f$ and $g$, as well as the characteristic functions of $X$ and $Y$).  We can obtain explicit formulas for these residuals by observing the optimality conditions for step 2 and 4 of Algorithm \ref{alg:pdhg}:
\begin{align}
0&\in F(x\kp)+A^Ty_k+\frac{1}{\tau_k}(x\kp-x_k)  \label{step1opt}\\
0&\in G(y\kp)-A(2x\kp-x_k)+\frac{1}{\sigma_k}(y\kp-y_k) . \label{step2opt}
\end{align}
Manipulating these optimality conditions yields
\begin{align}
P\kp &= \frac{1}{\tau_k}(x_k-x\kp) -A^T(y_k-y\kp) \\
	&\in\, F(x\kp)+A^Ty\kp \label{primal},\\
D\kp &= \frac{1}{\sigma_k}(y_k-y\kp)-A(x_k-x\kp) \\
	&\in \,G(y\kp)-A x\kp. \label{dual}
\end{align}

Formulas \eqref{primal} and \eqref{dual}  define a sequence of primal and dual residual vectors such that $P_k \in P(x_k,y_k)$ and $D_k \in D(x_k,y_k).$ 

We say that Algorithm \ref{alg:pdhg} converges if  
$$\lim_{k\to \infty} \|P_k\|^2 + \|D_k\|^2 = 0.$$
Note that studying the residual convergence of Algorithm \ref{alg:pdhg} is more general than studying the convergence of subsequences of iterates, because the solution to \eqref{saddle} need not be unique.

\section{Common Saddle-Point Problems} \label{section:problems}
Many common variational problems have saddle-point formulations that are efficiently solved using PDHG. While the applications of saddle-point problems are vast, we focus here on several simple problems that commonly appear in image processing.  
\subsection{Total-Variation Denoising}
A ubiquitous problem in image processing is minimizing the Rudin-Osher-Fatemi (ROF) denoising model \cite{ROF92}:
   \eqb{rof}\min_x | \nabla x| + \mot\|x-f \|^2.\eqe 
The energy \eqref{rof} contains two terms.  The   $\ell_2$ term on the right minimizes the squared error between the recovered image and the noise-contaminated measurements, $f$.  The TV term, $|\nabla x|,$ enforces that the recovered image be smooth in the sense that its gradient has sparse support. 

The problem\eqref{rof} can be formulated as a saddle-point problem by writing the TV term as a maximization over the ``dual'' variable $y\in \reals^{2\times N}$, where the image $x\in \reals^N$ has $N$ pixels:
 \eqb{tvmax}TV(x) = |\nabla x| = \max_{\|y\|_\infty\le 1}  y\cdot \nabla x.\eqe
Note that the maximization is taken over the set 
$$C_\infty =\{y\in R^{2\times N} \,| \,\, y_{1,i}^2+y_{2,i}^2\le1 \}.$$
The TV-regularized problem \eqref{rof} can then  be written in its primal-dual form
 \eqb{rof_pd}\max_{y\in  C_{\infty} } \min_x  \,  \mot \|x-f \|^2 +y\cdot \nabla x \eqe 
which is clearly a saddle-point problem in the form \eqref{saddle} with $f=\frac{\mu}{2}\|x-f \|^2$, $A=\nabla,$ $g =0,$ $X = \reals^N,$ and $Y = C_\infty$.

  To apply Algorithm \ref{alg:pdhg}, we need efficient solutions to the sub-problems in steps \ref{proxF} and \ref{proxG}, which can be written
    \begin{align}
    J_{\tau F}(\hat x) &=\argmin_x \mot\|x-f \|^2+\oott \|x- \hat x\|^2\\
    		 &= \frac{\tau}{\tau\mu+1}(\mu f + \tinv \hat x)   \\
    J_{\sigma G}(\hat y) &=\argmin_{y\in C_\infty} \oots \|y- \hat y\|^2 = \left ( \frac{y_i}{\max\{y_i,1\}} \right )_{i=1}^M. \label{shrink}
      \end{align}

\subsection{TVL1}
   Another common denoising model is TVL1 \cite{CE05}.  This model solves
   \eqb{tvl1}\min_x | \nabla x| + \mu |x-f |\eqe 
which is similar to \eqref{rof}, except that the data term relies on the $\ell_1$ norm instead of $\ell_2.$  This model is very effective for problems involving ``heavy-tailed'' noise, such as shot noise or salt-and-pepper noise. 

To put the energy \eqref{tvl1} into the form \eqref{saddle}, we write the data term in its variational form
 $$\mu |x-f| = \max_{y \in C_\mu }  y^T( x-f )$$ 
 where $C_\mu =\{y\in R^N \,| \,\, |y_{i}| \le \mu \}.$

 The energy \eqref{tvl1} can then be minimized  using the formulation
 \eqb{tvl1Saddle}
 \max_{y_1\in C_\infty,\, y_2\in  C_\mu}\min_x    y_1\cdot \nabla x + y_2^T x - y_2^T f  
 \eqe 
which is of the form \eqref{saddle}.
\subsection{Globally Convex Segmentation}
Image segmentation is the task of grouping pixels together by intensity and spatial location.   
One of the simplest models for image segmentation is the globally convex segmentation model of Chan, Esedoglu, and Nikolova (CEN) \cite{CEN06,GBO10}. Given an image $f$ and two real numbers $c_1$ and $c_2$ we wish to partition the pixels into two groups depending on whether their intensity lies closer to $c_1$ or $c_2.$  Simultaneously, we want the regions to have smooth boundaries.
  
 The CEN segmentation model has the variational form
\eqb{gcs}
\min_{0\le x\le 1}  \tv +  x^T  l
\eqe
where $l_i = (f_i-c_1)^2-(f_i-c_2)^2.$  The inner-product term of the right forces the entries in $x$ toward either 1 or 0, depending on whether the corresponding pixel intensity in $f$ is closer to $c_1$ or $c_2.$  At the same time, the TV term enforces that the boundary is smooth.  

Using the identity \eqref{tvmax}, we can write the model \eqref{gcs} as
\eqb{gcsSaddle}
\max_{y\in C_\infty} \min_{x\in [0,1]} \,  \tvmax +  x^T l.
\eqe
We can then apply Algorithm \ref{alg:pdhg}, where step \ref{proxF} takes the form
 $$J_{\tau F}(x) = \max\{0,\min\{x,1\}\}$$
  and step  \ref{proxG}  is given by \eqref{shrink}.

Generalizations of \eqref{gcs} to multi-label segmentations and general convex regularizers have been presented in \cite{GBO12, BCB12, GC10, BYT11}.  Many of these models result in minimizations similar to \eqref{gcs}, and PDHG methods have become a popular scheme for solving these models. 

\subsection{Compressed Sensing / Single-Pixel Cameras}
  In compressed sensing, one is interested in reconstructing images from incomplete measurements taken in an orthogonal transform domain (such as the Fourier or Hadamard domain).  It has been shown that high-resolution images can be obtained from incomplete data as long as the image has a sparse representation \cite{CRT06, CR05}, for example a sparse gradient.
  
    The Single Pixel Camera (SPC) is an imaging modality that leverages compressed sensing to reconstruct images using a small number of measurements \cite{DDTLSKB08}.  Rather than measuring in the pixel domain like conventional cameras, SPC's measure a subset of the Hadamard transform coefficients of an image.  
    
      To reconstruct images, SPC's rely on variational problems of the form
      \eqb{spc}
     \min_x \tv +\mot \|RH x - b \|^2
      \eqe 
where $b$ contains the measured transform coefficients of the image, and $H$ is an orthogonal transform matrix (such as the Hadamard transform). Here, $R$ is a diagonal ``row selector'' matrix, with diagonal elements equal to 1 if the corresponding Hadamard row has been measured, and 0 if it has not.  
 
 This problem can be put into the saddle-point form
 $$\max_{y \in C_\infty}\min_x    \mot \|RH x-b \|^2 +y\cdot \nabla x $$
 and then solved using PDHG.

The solution to step \ref{proxG} of Algorithm \ref{alg:pdhg} is given by \eqref{shrink}.  Step \ref{proxF} requires we solve 
 $$\min_x    \mot \|RH x-b \|^2  +\oott \|x- \hat x\|^2.$$
 Because $H$ is orthogonal, we can write the solution to this problem explicitly as
 $$x =  H^T\left(\mu R+\tinv I\right)^{-1}H\left(\mu H^TRb+\tinv\hat x\right).$$
 
 Note that for compressed sensing problems of the form \eqref{spc},  PDHG has the advantage that all steps of the method can be written explicitly, and no expensive minimization sub-steps are needed.    This approach to solving \eqref{spc} generalizes to any problem involving a sub-sampled orthogonal transform.  For example, the Hadamard transform $H$ could easily be replaced by a Fourier  transform matrix.

\subsection{$\ell_\infty$ Minimization}   
In wireless communications, one is often interested in finding signal representations with low peak-to-average power ratio (PAPR) \cite{HL05, SYB12}.  Given a signal $z$ with high dynamic range, we can obtain a more efficient representation by solving
  \eqb{linf} 
        \min_x\, \|x \|_\infty \text{ subject to } \| Dx - z\| \le \epsilon 
   \eqe
  where $x$ is the new representation, $D$ is the representation basis, and $\epsilon$ is the desired level of representation accuracy  \cite{SYB12}.  Minimizations of the form \eqref{linf} also arise in numerous other applications, including robotics \cite{Cadzow71} and nearest neighbor search \cite{JFF12}.

 The problem \eqref{linf} can be written in the constrained form
 $$\min_{x_1\in \reals^N,x_2\in C_\epsilon} \|x_1\|_\infty  \st x_2 = Dx_1-z  $$
 where $C_\epsilon = \{z \,| \,\, \|z\|<\epsilon \}.$  When we introduce Lagrange multipliers $y$ to enforce the equality constraint, we arrive at the saddle-point formulation
    $$\max_{y\in \reals^M} \min_{x_1\in \reals^N,x_2\in C_\epsilon} \|x_1\|_\infty  + y^T (Dx_1- x_2-z )  $$
    which is of the form \eqref{saddle}.  
    
    It often happens that $D$ is a complex-valued frame (such as a sub-sampled Fourier matrix).  In this case, $x$ and $y$ can take on complex values.  The appropriate saddle-point problem is then
          $$\max_{y\in \reals^M} \min_{x_1\in \reals^N,x_2\in C_\epsilon} \|x_1\|_\infty  + \Re\{  y^T ( Dx_1- x_2-z ) \}  $$
  and we can apply PDHG just as we did for real-valued problems provided we are careful to use the Hermitian definition of the matrix transpose (see the remark at the end of Section \ref{section:backtrack}).
  
     To apply PDHG to this problem, we need to evaluate the proximal minimization
     $$\min_x    \|x \|_\infty +\frac{1}{2t}\| x - \hat x\|^2$$
     which can be computed efficiently using a sorting algorithm \cite{DSSC08}.  
     
\subsection{Linear Programming}
 Suppose we wish to solve a large-scale linear program in the canonical form
\begin{align} \begin{split} \label{lp}
 &\min_x \quad c^T  x  \\
 &\text{ subject to\quad} Ax\le b, \, x\ge0
 \end{split}\end{align}
where  $c$ and $b$ are vectors, and $A$ is a matrix \cite{Schrijver03}.  This problem can be written in the primal-dual form
$$\max_{y\ge 0} \min_{x\ge 0}   c^T x  +  y^T( Ax-b )  $$
where $y$ can be interpreted as Lagrange multipliers enforcing the condition $Ax\le b$ \cite{Schrijver03}.

Steps \ref{proxF} and \ref{proxG} of Algorithm \ref{alg:pdhg} are now simply
     \begin{align*}
    J_{\tau F}(\hat x) &=  \max\{ \hat x - \tau c ,0\},\\
      J_{\sigma G}(\hat y) &=\max\{ \hat y - \sigma b ,0\}.
            \end{align*}

The PDHG solution to linear programs is advantageous for some extremely large problems for which conventional simplex and interior-point solvers are intractable. Primal-dual solvers for \eqref{lp} have been studied in \cite{PC11}, and compared to conventional solvers.  It has been observed that this approach is not efficient for general linear programming, however it behaves quite well for certain structured programs, such as those that arise in image processing \cite{PC11}.

\section{Residual Balancing}
  When choosing the stepsize for PDHG, there is a tradeoff between the primal and dual residuals.   Choosing a large value of $\tau_k$ creates a very powerful minimization step in the primal variables and a slow maximization step in the dual variables, resulting in very small primal residuals at the cost of large dual residuals.  Choosing $\tau_k$ to be small, on the other hand, results in small dual residuals at the cost of large primal errors. 
  
   Ideally, one would like to choose stepsizes so that the larger of $P_k$ and $D_k$ is as small as possible.  If we assume the primal/dual residuals decrease/increase monotonically with $\tau_k,$ then $\max\{P_k,D_k\}$ is minimized when both residuals are equal in magnitude.    This suggests that $\tau_k$ be chosen to ``balance'' the primal and dual residual -- i.e., the primal and dual residuals should be roughly the same size, up to some scaling to account for units.  This principle has been suggested for other iterative methods (see \cite{HYW00,BPCPE10}).  
   
   The particular residual balancing principle we suggest is to enforce
    \eqb{balance}|P_k| \approx s|D_k|  \eqe 
    where $|\cdot|$ denotes the $\ell_1$ norm, and $s$ is a scaling parameter.  Residual balancing methods work by ``tuning'' parameters after each iteration to approximately maintain this equality.  If the primal residual grows disproportionately large, $\tau_k$ is increased and $\sigma_k$ is decreased (or vice-versa).  
     
  In \eqref{balance} we use the $\ell_1$-norm to measure the size of the residuals.  In principle, $\ell_1$ could easily be replaced with the $\ell_2$ norm, or any other norm.  We have found that the $\ell_1$ norm performs somewhat better than the $\ell_2$ norm for some problems, because it is less sensitive to outliers that may dominate the residual. 

  Note the proportionality in \eqref{balance} depends on a constant $s$.  This is to account for the effect of scaling the inputs (i.e., changing units).  For example, in the TVL1 model \eqref{tvl1} the input image $f$ could be scaled with pixel intensities in $[0,1]$ or $[0,255].$  As long as the primal step sizes used with the latter scaling are $255$ times that of the former (and the dual step sizes are 255 times smaller) both problems produce similar sequences of iterates. However, in the $[0,255]$ case the dual residuals are 255 times larger than in the $[0,1]$ case while the primal residuals are the same.  

  For image processing problems, we recommend using the $[0,255]$ scaling with $s=1,$ or equivalently the $[0,1]$ scaling with $s=255.$ This scaling enforces that the saddle point objective \eqref{saddle} is nearly maximized with respect to $y$ and that the saddle point term \eqref{tvmax} is a good approximation of the total variation semi-norm.

\section{Adaptive Methods}
In this section, we develop new adaptive PDHG methods.  The first method assumes we have a bound on $\rho(A^TA).$  In this case, we can enforce the stability condition 
\eqb{stab}\tau_k \sigma_k <L<\rho(A^TA)^{-1}.
\eqe
  This is the same stability condition that guarantees convergence in the non-adaptive case.  In the adaptive case, this condition alone is not sufficient for convergence but leads to relatively simple methods.  

  The second method does not require any knowledge of $\rho(A^TA).$  Rather, we use a backtracking scheme to gaurantee convergence.  This is similar to the Armijo-type backtracking line searches that are used in conventional convex optimization.  In addition to requiring no knowledge of $\rho(A^TA),$ the backtracking scheme has the advantage that it can use relatively long steps that violate the stability condition \eqref{stab}.   Especially for problems with highly sparse solutions, this can lead to faster convergence than methods that enforce \eqref{stab}.

 Both of these methods have strong convergence guarantees.  In particular, we show in Section \ref{section:theory} that both methods converge in the sense that the norm of the residuals goes to zero.
 
 \subsection{Adaptive PDHG}
 
   The first adaptive method is listed in Algorithm \ref{alg:adapt}.
\begin{algorithm*}
    \caption{Adaptive PDHG}
    \label{alg:adapt}
    \begin{algorithmic}[1]
    \Require $x_0\in \reals^N, y_0\in \reals^M, \sigma_0\tau_0<\rho(A^TA)^{-1} , (\alpha_0,\eta) \in (0,1)^2, \Delta>1,s>0$ 
    \While{$p_k,d_k>{tolerance}$  }
    \State $x\kp = J_{\tau_k F} ( x_k - \tau_k A^Ty_k) $ \Comment{Begin with normal PDHG}
     \State $y\kp = J_{\sigma_k G} ( y_k + \sigma_k A (2 x\kp-x_k) ) $  
      \State  $p\kp = |(x_k-x\kp)/\tau_k - A^T(y_k-y\kp)|$ \Comment{Compute primal residual}\label{adapt:primal}
   \State $d\kp =  |(y_k-y\kp)/\sigma_k - A(x_k-x\kp)|$  \Comment{Compute dual residual}\label{adapt:dual}
        
    \If {$p\kp> sd\kp \Delta $}       \Comment{If primal residual is large...}
		\State    $ \tau\kp = \tau_k/(1-\alpha_k)$  \Comment{ Increase primal stepsize}
		\State    $ \sigma\kp = \sigma_k(1-\alpha_k)$  \Comment{Decrease dual stepsize}
		\State    $ \alpha\kp = \alpha_k \eta $          \Comment{Decrease adaptivity level}
    
    \EndIf
  \If {$p\kp<sd\kp/\Delta$}    \Comment{If dual residual is large...}
		\State    $ \tau\kp = \tau_k(1-\alpha_k)$    \Comment{Decrease primal stepsize}
		\State    $ \sigma\kp = \sigma_k/(1-\alpha_k)$ \Comment{Increase dual stepsize}
		\State    $ \alpha\kp = \alpha_k \eta $    \Comment{Decrease adaptivity level}
    \EndIf
 \If {$  sd\kp /\Delta\le p\kp\le sd\kp \Delta$}    \Comment{If residuals are similar...}
		\State    $ \tau\kp = \tau_k$    \Comment{Leave primal step the same}
		\State    $ \sigma\kp = \sigma_k$ \Comment{Leave dual step the same}
		\State    $ \alpha\kp = \alpha_k  $    \Comment{Leave adaptivity level the same}
    \EndIf
    \EndWhile
       \end{algorithmic} 
  \end{algorithm*}
The loop in Algorithm \ref{alg:adapt} begins by performing a standard PDHG step using the current stepsizes.    In steps \ref{adapt:primal} and \ref{adapt:dual}, we compute the primal and dual residuals and store their $\ell_1$ norms in $p_k$ and $d_k.$   If the primal residual is sufficiently large compared to the dual residual then the primal stepsize is increased by a factor of $(1-\alpha_k)^{-1},$ and the dual stepsize is decreased by a factor of $(1-\alpha_k).$  If the primal residual is somewhat smaller than the dual residual then the primal stepsize is decreased and the dual stepsize is increased.  If both residuals are comparable in size, then the stepsizes remain the same on the next iteration.  

   The parameter $\Delta>1$ is used to compare the sizes of the primal and dual residuals.  The stepsizes are only updated if the residuals differ by a factor greater than $\Delta.$
   
   The sequence $\{\alpha_k\}$ controls the adaptivity level of the method.  We start with some $\alpha_0\in (0,1).$  Every time we choose to update the stepsize parameters, we define $\alpha\kp = \eta \alpha_k$ for some $\eta<1.$  In this way, the adaptivity decreases over time.  We will show in Section \ref{section:theory} that this definition of $\{\alpha_k\}$ is needed to guarantee convergence of the method.

 \begin{remark}
 The adaptive scheme requires several arbitrary constants as inputs.  These are fairly easy to choose in practice.  A fairly robust choice is $a_0 = 0.5,$ $\Delta = 1.5,$ $\eta = 0.95,$ and $\tau_0=\sigma_0 = 1/\sqrt{\rho(A^TA)}.$  However, these parameters can certainly be tuned for various applications.  As stated above, we recommend choosing $s=1$ for imaging problems when pixels lie in $[0,255]$, and $s=255$ when pixels lie in $[0,1].$
 \end{remark}
 
  \begin{remark}
 The computation of the residuals in steps \ref{adapt:primal} and \ref{adapt:dual} of Algorithm \ref{alg:adapt} requires multiplications by $A^T$ and $A$, which seems at first to increase the cost of each iteration.  Note however that the values of $A^Ty_k$ and $Ax_k$ are already computed in the other steps of the algorithm.  If we simply note that $A^T(y_k-y\kp) = A^Ty_k-A^Ty\kp$ and $A(x_k-x\kp) = Ax_k-Ax\kp,$ then we can evaluate the residuals in steps   \ref{adapt:primal} and \ref{adapt:dual} without any additional multiplications by $A$ or $A^T$.  
 \end{remark}

 \subsection{Backtracking PDHG} \label{section:backtrack}
 
   A simple modification of Algorithm \ref{alg:adapt} allows the method to be applied when bounds on $\rho(A^TA)$ are unavailable,  or when the stability condition \eqref{stab} is overly conservative.  This is accomplished by choosing ``large'' initial stepsizes,  and then decreasing the stepsizes using the ``backtracking'' scheme described below.

We will see in Section \ref{section:theory} that convergence is guaranteed if the following ``backtracking'' condition holds at each step:
\eqb{back}
b_k = \frac{2\tau_k\sigma_k (y\kp-y_k )^T A(x\kp-x_k)} {\gamma\sigma_k\|x\kp-x_k\|^2+\gamma\tau_k\|y\kp-y_k\|^2}\le 1.
\eqe
where $\gamma\in (0,1)$ is a constant. If this inequality does \emph{not} hold, then the stepsizes are too large.  In this case we decrease the stepsizes by  choosing 
\eqb{backupdates} 
\tau\kp = \beta\,\tau_k/b_k, \text{ and }  \sigma\kp =  \beta\,\sigma_k/b_k,
\eqe
 for some $\beta \in (0,1).$ We recommend choosing $\gamma=0.75$ and $\beta = 0.95$, although the method is convergent for any $\gamma,\beta\in (0,1).$   Reasons for this particular update choice, are elaborated in the Section \ref{section:backTheory}.

We will show in Section \ref{section:theory} that this backtracking step can only be activated a finite number of times.   Furthermore, enforcing the condition \eqref{back} at each step is sufficient to guarantee convergence, provided that one of the spaces $X$ or $Y$ is bounded.

  The initial stepsizes $\tau_0$ and $\sigma_0$ should be chosen so that $\tau_0\sigma_0$ is somewhat larger than $\rho(A^TA)^{-1}.$   In the numerical experiments below, we choose
  \eqb{init}
  	\tau_0=\sigma_0 = \sqrt{\frac{2\|x_r\|}{\|A^TAx_r\|}}
  \eqe
 where $x_r$ is a random Gaussian distributed vector.  Note that the ratio $\frac{\|A^TAx_r\|}{\|x_r\|}$  is guaranteed to be less than $\rho(A^TA),$ and thus  $\tau_0\sigma_0>\rho(A^TA)^{-1}.$  The factor of 2 is included to account for problems for which the method is stable with large steps.

\begin{remark}
 Note that our convergence theorems require either $X$ or $Y$ to be bounded.  This requirement is indeed satisfied by most of the common saddle-point problems described in Section \ref{section:problems}.  However, we have found empirically that the scheme is quite stable even for unbounded $X$ and $Y$. For the linear programming example, neither space is bounded.  Nonetheless, the backtracking method converges reliably in practice. 
  \end{remark}

\begin{remark}
In the case of complex-valued problems, both $x$ and $y$ may take on complex values, and the relevant saddle-point problem is
$$\min_{x\in X} \max_{y\in Y} f(x)+\Re\{ y^T Ax \} - g(y)$$
where $\Re\{\cdot\}$ denotes the real part of a vector.
In this case the numerator of $b_k$ may have an imaginary component. The algorithm is still convergent as long as we replace $b_k$ with its real part and use the definition
\eqb{backComp}
b_k = \frac{\Re\left \{2\tau_k\sigma_k (y\kp-y_k )^T A(x\kp-x_k)  \right \} } {\gamma\sigma_k\|x\kp-x_k\|^2+\gamma\tau_k\|y\kp-y_k\|^2}.
\eqe
  \end{remark}
  
\section{Convergence Theory} \label{section:theory}
 \subsection{Variational Inequality Formulation}
 
 For notational simplicity, we define the vector quantities 
 \eqb{vecs}
  u_k = \matb
  x_k\\
    y_k 
\mate,
\qquad
 R(u) = \matb
   P(x,y)\\
   D(x,y) \\
\mate,
\eqe
and the matrices
\begin{align}
\begin{split}\label{operators}
\qquad
 M_k = \matb
   \tau^{-1}_kI & -A^T\\
    -A & \sigma^{-1}_k I \\
\mate, \,\,\,
 &H_k = \left(\begin{matrix}
   \tau^{-1}_k I & 0\\
    0 & \sigma^{-1}_k I\\
\end{matrix}\right),
\\
Q(u)& = \left(\begin{matrix}
  A^Ty\\
   -Ax \\
\end{matrix}\right).
\end{split}\end{align}

This notation was first suggested to simplify PDHG by He and Yuan \cite{HY12}.

Using this notation, it can be seen that the iterates of PDHG satisfy
\eqb{vi}
0 \in R(u\kp) + M_k(u\kp-u_k) .
\eqe
Also, the optimality conditions \eqref{opt1} and \eqref{opt2} can be written succinctly as
\eqb{viOpt}
0\in R(u\opt) .
\eqe
Note that $R$ is monotone (see \cite{BC11}), meaning that
$$ ( u - \hat u )^T ( R(u) - R(\hat u) ) \ge  0, \quad \forall u,\hat u.$$  

Following \cite{HY12}, it is also possible to formulate each step of PDHG as a variational inequality (VI).
If  $u\opt=(x\opt,y\opt)$ is a solution to \eqref{saddle}, then $x\opt$  minimizes \eqref{saddle} (for fixed $y\opt$) .  More formally,
\eqb{VI1}
f(x)-f(x\opt)+ ( x-x\opt)^T A^T y\opt  \ge \quad \forall \, x \in X.
\eqe
Likewise, \eqref{saddle} is maximized by $y\opt,$ and so
\eqb{VI2}
-g(y)+g(y\opt)+ (y-y\opt)^T A x\opt  \le 0 \quad \forall \, y \in Y.
\eqe
Subtracting \eqref{VI2} from \eqref{VI1} yields the VI formulation
\eqb{VI}
\phi(u)-\phi( u\opt) + (u-u\opt)^TQ(u\opt) \ge 0 \quad\forall u\in \Omega. 
\eqe
where $\phi(u) = f(x)+g(y)$ and $\Omega = X\times Y.$
We say that a point $\tilde u$ is an approximate solution to \eqref{saddle} with VI accuracy $\epsilon$ if
$$\phi(u)-\phi( \tilde u) + (u-\tilde u)^TQ(\tilde u) \ge -\epsilon \quad\forall u\in B_1(\tilde u) \cap \Omega $$
where $B_1(\tilde u)$ is a unit ball centered at $\tilde u.$

In Section \ref{sec:thrms}, we prove two convergence results for adaptive PDHG.  First we prove the residuals vanish as $k\to \infty$.  Next, we prove a $O(1/k)$ ergodic convergence rate using the VI notion of convergence.

\subsection{Stepsize Conditions} \label{sec:cond}
We now discuss conditions that guarantee convergence of adaptive PDHG.  We delay the formal convergence proofs until Section \ref{section:proofs}.   We begin by defining the following constants:
\begin{align}
\delta_k &= \min\left \{\frac{\tau\kp}{\tau_k},\frac{\sigma\kp}{\sigma_k},1 \right \} \nonumber\\
\phi_k &= 1-\delta_k= \max \left\{ \frac{\tau_k-\tau\kp}{\tau_k}, \frac{\sigma_k-\sigma\kp}{\sigma_k},0 \right\}. \label{def:phi}
\end{align}
The constant $\delta$ quantifies the stepsize change, and $\phi$ quantifies the relative decrease in the stepsizes between iterations.

 Consider the following conditions on the step-sizes in Algorithm \ref{alg:pdhg}.

\begin{center}
\begin{figure*}{\linewidth-1cm}	
 \begin{framed}
\noindent \textbf{Convergence Conditions for Adaptive PDHG}

Algorithm \ref{alg:pdhg} converges if the following three requirements hold:
\begin{description}
\item [A] The sequences $\{\tau_k\}$ and $\{\sigma_k\}$ are bounded.
\item[B]  The sequence $\{\phi_k\}$ is summable, i.e., for some $C_\phi$
$$\sum_{k\ge0} \phi_k<C_\phi<\infty. $$
\item[C] One of the following two conditions is met:
\begin{description}
\item[C1] There is a constant $L$ such that for all $k>0$
$$\tau_k\sigma_k  <L< \rho(A^TA)^{-1}.  $$
\item[C2] Either $X$ or $Y$ is bounded, and there is a \\constant $c\in (0,1)$ such that for all $k>0$
$$\| u\kp-u_k \|^2_{M_k}\ge c\|u\kp-u_k\|^2_{H_k}.  $$
\end{description}
\end{description}
 \end{framed}
\end{figure*}
\end{center}

Two conditions must always hold to guarantee convergence:  (\A{})  The stepsizes must remain bounded, and (\B{}) the sequence $\{\phi_i\}$ must be summable.  Together, these conditions ensure that the steps sizes do not oscillate too wildly as $k$ gets large.  

  In addition,  either \Cone{} or \Ctwo{} must hold.   When we know the spectral radius of $A^TA,$ we can use condition \Cone{} to enforce that the method is stable.  When we have no knowledge of $\rho(A^TA),$ or when we want to take large stepsizes, the backtracking condition \Ctwo{} can also be used to enforce stability.  We will see that for small enough stepsizes, this backtracking condition can always be satisfied.
 
 Note that condition \Cone{} is slightly stronger than condition \Ctwo{}.  The advantage of condition \Cone{} is that it results is somewhat simpler methods because backtracking does not need to be used.  However when backtracking is used to enforce condition \Ctwo{}, we can take larger steps that violate condition \Cone{}.

In the following sub-sections, we explain how Algorithm \ref{alg:adapt} and its backtracking variant explicitly satisfy conditions \A{}, \B{}, and \C{}, and are therefore guaranteed to converge.  In Section \ref{section:proofs}, we prove convergence of the general Algorithm \ref{alg:pdhg} under assumptions  \A{}, \B{}, and \C{}.   
   
 \subsection{Convergence of the Adaptive Method }  \label{section:proofs}
   Algorithm \ref{alg:adapt} is a practical implementation of adaptive PDHG that satisfies conditions \A{}, \B{}, and \Cone{}.   
   An examination of Algorithm \ref{alg:adapt} reveals that there are three possible stepsize updates depending on the balance between the residuals.  Regardless of which update occurs, the product $\tau_k\sigma_k$ remains unchanged at each iteration, and so $\tau_k\sigma_k=\tau_0\sigma_0=L<\rho(A^TA)^{-1}.$ It follows that condition \A{} and \Cone{} are satisfied.  
   
   Also, because $0<\alpha_k<1$, we have 
  
   \begin{align}\begin{split}
   \phi_k &=1-\min \left \{\frac{\tau\kp}{\tau_k},\frac{\sigma\kp}{\sigma_k},1\right \}\\
   &=\begin{cases}
    \alpha_k \text{,  if the steps are updated}\\
     1 \text{,  otherwise.}
   \end{cases}
   \end{split}\end{align}
   Note that every time we change the stepsizes, we update $\alpha\kp = \eta \alpha_k$ where $\eta<1.$ Thus the non-zero entries in the sequence $\{\phi_k\}$ form a decreasing geometric sequence.  It follows that the summation condition \B{} is satisfied. 
  
    Because Algorithm \ref{alg:adapt} satisfies conditions  \A{} , \B{} and \C{}, its convergence is guaranteed by the theory presented in Section \ref{section:proofs}.  
   
 \subsection{Convergence of the Backtracking Method} \label{section:backTheory}
  As we shall prove in Section \ref{section:proofs}, the backtracking step is triggered only a finite number of times.  Thus, the backtracking method satisfies conditions \A{} and \B{}.  
   
   We can expand condition \Ctwo{} using the definition of $\|\cdot \|_{M_k}$ and  $\|\cdot \|_{H_k}$ to obtain
\begin{align*}
\frac{1}{\tau_k}\|x\kp-x_k\|^2 +&\frac{1}{\sigma_k}\|y\kp-y_k\|^2\\
		&-  2  (y\kp-y_k )^T A(x\kp-x_k)  \\
  >&  \frac{c}{\tau_k}\|x\kp-x_k\|^2 +\frac{c}{\sigma_k}\|y\kp-y_k\|^2. 
\end{align*}
If we combine like terms and let $\gamma = 1-c,$ then this condition is equivalent to
\begin{align} \label{back2}
\frac{\gamma}{\tau_k}\|x\kp-&x_k\|^2  +\frac{\gamma}{\sigma_k}\|y\kp-y_k\|^2\\
 &>  2 (y\kp-y_k )^T A(x\kp-x_k) .
\end{align}
If we note that the left side of \eqref{back2} is non-negative, then we can easily see that this condition is equivalent to \eqref{back}. Thus, the backtracking conditions \eqref{back} explicitly enforce condition \Ctwo{}, and the convergence of Algorithm \ref{alg:adapt} with backtracking is guaranteed by the analysis below. 

We now address the form of the backtracking update \eqref{backupdates}.
A simple Armijo-type backtracking scheme would simply choose $\tau\kp = \xi_k \tau_k$ and $\sigma\kp = \xi \sigma_k,$ where $\xi<0,$ every time that \eqref{back2} is violated.  However if our initial guess for $L$ is very large, then this backtracking scheme could be very slow.  Rather, we wish to predict the value of $\xi$ that makes \eqref{back2} hold.  To do this, we assume (as a hueristic) that the values of $\|x\kp-x_k\|$  and $\|y\kp-y_k\|$ decrease linearly with $\xi.$  Under this assumption, the value of $\xi$ that makes \eqref{back2} into an equality is $\xi_k=b_k^{-1}.$  In order to guarantee that the stepsizes can become arbitrarily small, we make the slightly more conservative choice of $\xi_k = \beta / b_k$ for $\beta<1.$

   \subsection{Sketch of Proof}

The overall strategy of our convergence proof is to show that the PDHG method satisfies a \emph{Fej\'er inequality}, i.e., an inequality stating that the iterates move toward the solution set on every iteration.  We derive a simple Fej\'er inequality in Lemma \ref{lemma:ineq}.  This lemma states that $u\kp$ always lies closer to the solution set than $u_k$ as measured in the $M_k$ norm.  Normally, this would be enough to prove convergence (see \cite{BC11}, chapter 5).  However, in our case the proof is not so straight-forward for several reasons.  First, the $M_k$-norm changes every time we update the stepsizes.  Second, when the backtracking condition \Ctwo{} is used, $M_k$ may be indefinite, in which case $\|\cdot\|_{M_k}$ is not a proper norm, and can even assume negative values.

  Nonetheless, we can still prove convergence.   We begin by proving that \Cone{} guarantees that $M_k$ is positive definite.  We then derive a Fej\'er inequality in Lemma \ref{lemma:ineq}.   In cases when $M_k$ is indefinite, Lemmas \ref{lemma:upper} and \ref{lemma:lower} show that terms involving the $M_k$-norm remain uniformly bounded from both above \emph{and} below.  These bounding conditions will be strong enough to allow us to complete the proof  in the case that $M_k$ is indefinite.

   \subsection{Preliminary Lemmas}
   We now prove several preliminary results about the iterates of PDHG that will be needed later for the convergence arguments. 
  This first lemma has been adapted from He and Yuan \cite{HY12}.  It shows that the stability condition \Cone{} forces the matrix $M_k$ to be positive definite.  When this condition is satisfied, the operator $ u^T M_ku =\|u\|_{M_k}^2$ is a proper norm and can be used in our convergence analysis.

 Lemma \ref{lemma:posdef} also shows that backtracking can only occur a finite number of times.  When the product $\tau_k\sigma_k$ becomes sufficiently small, so does the constant $C_M$ in the lemma, and the ratio \eqref{back} will always be less than 1.

\begin{lemma} \label{lemma:posdef}
If $\tau\sigma\|A^TA\|_{op}<1$ then $M$ is positive definite, and we have
\begin{align*}
\left(1-C_M\right)\left(  \frac{1}{\sigma}\|y\|^2 +\frac{1}{\tau}\|x\|^2   \right)\le \|u\|_M^2 \\
  \le \left(1+C_M\right)\left(  \frac{1}{\sigma}\|y\|^2 +\frac{1}{\tau}\|x\|^2   \right)
  \end{align*}
where $C_M = \sqrt{\tau\sigma\|A^TA\|_{op}}.$
\end{lemma}
\begin{IEEEproof}
By definition
\begin{align}\label{posdef:ineq1}
 u^T Mu &=\frac{1}{\tau}\|x\|^2 - 2 y^TAx  + \frac{1}{\sigma}\|y\|^2. 
\end{align}
Now use the inequality $a^2+b^2>2ab$ to obtain
\begin{align} 
\begin{split}\label{posdef:ineq2}
 2  y^T Ax &\le 2\frac{\|A\|_{op}\sqrt\tau\|y\|}{(\tau\sigma\|A^TA\|_{op})^{\frac{1}{4}}}\frac{(\tau\sigma\|A^TA\|_{op})^{\frac{1}{4}}}{\sqrt\tau}\|x\| \\
 &\le  \frac{\|A^TA\|_{op} \tau}{\sqrt{\tau\sigma\|A^TA\|_{op}}} \|y\|^2 + \frac{\sqrt{\tau\sigma\|A^TA\|_{op}}}{\tau}\|x\|^2\\
 &= \sqrt{\tau\sigma\|A^TA\|_{op}}\left(  \frac{1}{\sigma}\|y\|^2 +\frac{1}{\tau}\|x\|^2   \right).
\end{split}
\end{align}
Applying \ref{posdef:ineq2} to \ref{posdef:ineq1} yeilds
$$\|u\|_M^2 \ge \left(1-\sqrt{\tau\sigma\|A^TA\|_{op}}\right)\left(  \frac{1}{\sigma}\|y\|^2 +\frac{1}{\tau}\|x\|^2   \right).$$
We also have that
\begin{align*}\begin{split} 
 u^T Mu =&\frac{1}{\tau}\|x\|^2 - 2 y^T Ax + \frac{1}{\sigma}\|y\|^2 \\
\le& \frac{1}{\tau}\|x\|^2 + 2\|A\|_{op}\|x\|\|y\|+ \frac{1}{\sigma}\|y\|^2\\
\le& \frac{1}{\tau}\|x\|^2+ \frac{1}{\sigma}\|y\|^2\\
&+\sqrt{\tau\sigma\|A^TA\|_{op}}\left(  \frac{1}{\sigma}\|y\|^2 +\frac{1}{\tau}\|x\|^2   \right)\\
\end{split}\end{align*}
and so
$$\|u\|_M^2 \le \left(1+\sqrt{\tau\sigma\|A^TA\|_{op}}\right)\left(  \frac{1}{\sigma}\|y\|^2 +\frac{1}{\tau}\|x\|^2   \right).$$
\end{IEEEproof}

We next show that the distance between the true solution and the PDHG iterates is decreasing in the $M_k$ norm.  If $M_k$ were constant and positive definite, then this result would be sufficient for convergence.  However, because the matrix $M_k$ changes at each iteration and may be indefinite, this condition does not guarantee convergence on its own.    
\begin{lemma} \label{lemma:ineq}
The iterates generated by Algorithm \ref{alg:pdhg} satisfy
\begin{align*}
\|u_k - u\opt\|^2_{M_k} &\ge \| u\kp - u_k\|^2_{M_k}+\| u\kp - u\opt\|^2_{M_k}.
\end{align*}
\end{lemma}
\begin{IEEEproof}
Subtracting \eqref{vi}  from \eqref{viOpt} gives us
$$
   M_k(u\kp-u_k) \in R(u\opt) -  R(u\kp) .
$$
Taking the inner product with $(u\opt - u\kp)$ gives us
\begin{align*}
( u\opt - u\kp)^T M_k&(u\kp-u_k)  \\
&\ge ( u\opt - u\kp)^T (R(u\opt) -  R(u\kp) ).
\end{align*}
Because $R$ is monotone, the right hand side of the above equation is non-negative, and so
\eqb{converge:ineq1}
(u\opt - u\kp)^T M_k(u\kp-u_k) \ge 0.
\eqe
Now, observe the identity
\begin{align*}
\|u_k - u\opt\|^2_{M_k} &= \| u\kp - u_k\|^2_{M_k}+\| u\kp - u\opt\|^2_{M_k}
\\&\hspace{1cm}+2 ( u_k-u\kp)^T  M_k( u\kp - u\opt) .
\end{align*}
Applying \eqref{converge:ineq1} yields the result.
\end{IEEEproof}

We now show that the iterates generated by PDHG remain bounded.

\begin{lemma} \label{lemma:upper}
Suppose  the step sizes for Algorithm \ref{alg:pdhg} satisfy conditions \A{}, \B{} and \C{}.
Then 
$$\|u_k-u\opt\|^2_{H_k}\le C_U$$
for some upper bound $C_U>0.$
\end{lemma}

\begin{IEEEproof}   We first consider the case of condition \Cone{}.
From \eqref{posdef:ineq2} we have
\begin{align*} \begin{split}
 2 y^T Ax &\le \sqrt{\tau\sigma\|A^TA\|_{op}}\left(  \frac{1}{\sigma}\|y\|^2 +\frac{1}{\tau}\|x\|^2  \right).
\end{split}\end{align*}
Subtracting $2 \sqrt{\tau\sigma\|A^TA\|_{op}}  y^TAx$ from both sides yields
\begin{align*} 
& \left(2-2\sqrt{\tau\sigma\|A^TA\|_{op}}\right) y^TAx\\
 &\hspace{40pt}\le \sqrt{\tau\sigma\|A^TA\|_{op}}\left(  \frac{1}{\sigma}\|y\|^2 +\frac{1}{\tau}\|x\|^2   -2 y^TAx\right)\\
 &\hspace{40pt}=\|u\|^2_M\sqrt{\tau\sigma\|A^TA\|_{op}}.
\end{align*}
Taking $x = x\kp-x\opt,$  $y = y\kp-y\opt,$ $\tau= \tau\kp,$ and $\sigma = \sigma\kp,$ and noting that $\tau\kp\sigma\kp<L$ we obtain 
$$
 2 (y\kp-y\opt)^T A(x\kp-x\opt)   \le  C_1\|u\kp-u\opt\|^2_{M\kp}
$$
where $C_1 = \frac{ \sqrt{L\|A^TA\|_{op}}}{1- \sqrt{L\|A^TA\|_{op}}}>0.$

Applying this to the result of Lemma \ref{lemma:ineq}, yields
\begin{align*} 
\|u_k - u&\opt\|^2_{M_k} \ge  \| u\kp - u\opt\|^2_{M_k}\\
=& \|u\kp-u\opt\|^2_{H_k}+ (y\kp-y\opt )^T A(x\kp-x\opt)\\
\ge& \delta_k\|u\kp-u\opt\|^2_{H\kp}+ (y\kp-y\opt)^T A(x\kp-x\opt)\\
\ge& \delta_k\|u\kp-u\opt\|^2_{M\kp}\\
&+(1-\delta_k) (y\kp-y\opt)^T A(x\kp-x\opt)\\
\ge& \delta_k\|u\kp-u\opt\|^2_{M\kp}-(1-\delta_k)C_1\|u\kp-u\opt\|^2_{M\kp}\\
=& ( 1-(1+C_1) \phi_k ) \|u\kp-u\opt\|^2_{M\kp}. 
\end{align*}
Note that $(1+C_1) \phi_k >1 $ for only finitely many $k$ , and so we assume without loss of generality that $k$ is sufficiently large that $( 1-(1+C_1) \phi_k )>0.$
We can then write 
\eqb{bound:ineq2}
\|u_1 - u\opt\|^2_{M_1}\ge \prod_{k=1}^{n-1} ( 1-(1+C_1) \phi_k ) \|u_n-u\opt\|^2_{M_n}.
\eqe
Since $\sum_k \phi_k<\infty$, we have that  $\sum_k (1+C_1)\phi_k<\infty$, and the product on the right of \eqref{bound:ineq2} is bounded away from zero.  Thus, there is a constant $C_2$ with
$$
\|u_1 - u\opt\|^2_{M_1}\|\ge C_2\|u_n-u\opt\|^2_{M_n} \ge C_2(1-C_M)\|u_n-u\opt\|^2_{H_n}
$$
and the lemma is true in the case of assumption \Cone{}.

 We now consider the case where condition \Ctwo{} holds.  We assume without loss of generality that $Y$ is bounded (the case of bounded $X$ follows by nearly identical arguments).  In this case, we have $\|y\|\le C_y$ for all $y\in Y.$

Note that
  \begin{align} \|u\kp-&u\opt \|_{M\kp} = - (y\kp-y\opt )^T A(x\kp-x\opt)  \\
  &+ \frac{1}{\tau\kp}\|x\kp-x\opt\|^2 +\frac{1}{\sigma\kp}\|y\kp-y\opt\|^2 \nonumber \\
  \ge& -2C_y\|A\|_{op} \|x\kp-x\opt\|   \\
  &+ \frac{1}{\tau\kp}\|x\kp-x\opt\|^2 + \frac{1}{\sigma\kp}\|y\kp-y\opt\|^2. \label{bound:upper} 
  \end{align}
  When $\|x\kp-x\opt\|$ grows sufficiently large, the term involving the square of this norm dominates the value of \eqref{bound:upper}.   Since $\{\tau_k\}$ and $\{\sigma_k\}$ are bounded from above, it follows that there is some positive $C_x$ such that whenever
    \begin{align}  \label{upper:cond}
     \frac{1}{\tau\kp}\|x\kp-x\opt\|^2 + \frac{1}{\sigma\kp}\|y\kp-y\opt\|^2 \ge C_x
       \end{align} 
  we have
    \begin{align} 
        \frac{1}{\tau\kp}\|x\kp-x\opt\|^2 &+ \frac{1}{\sigma\kp}\|y\kp-y\opt\|^2 \\
             & \ge 4C_y\|A\|_{op} \|x\kp-x\opt\| . \label{bound:1}  
        \end{align} 
        Combining \eqref{bound:1} with \eqref{bound:upper} yields
         \begin{align} 
             2\|u\kp-u\opt \|_{M\kp}&\ge  \frac{1}{\tau\kp}\|x\kp-x\opt\|^2\\
            & + \frac{1}{\sigma\kp}\|y\kp-y\opt\|^2 
   \end{align}
   whenever \eqref{upper:cond} holds.  In this case, we have
  \begin{align}
   \|u\kp-&u\opt \|^2_{M_k} =- (y\kp-y_k )^T A(x\kp-x\opt) \\
     &+ \frac{1}{\tau_k}\|x\kp-x\opt\|^2 + \frac{1}{\sigma_k}\|y\kp-y\opt\|^2 \nonumber \\
    \ge&  - (y\kp-y\opt)^T A(x\kp-x\opt) \\
      &+  \frac{\delta_k}{\tau\kp}\|x\kp-x\opt\|^2+ \frac{\delta_k}{\sigma\kp}\|y\kp-y\opt\|^2 \nonumber\\
      =& \|u\kp-u\opt \|^2_{M\kp}\\ 
       &-  \frac{\phi_k }{\tau\kp}\|x\kp-x\opt\|^2 - \frac{\phi_k }{\sigma\kp}\|y\kp-y\opt\|^2 \nonumber \\    
         &\ge (1-2\phi_k) \|u\kp-u\opt \|^2_{M\kp}.  \label{biggerThanC}
    \end{align}
    Applying \eqref{biggerThanC} to Lemma \ref{lemma:ineq}, we see that
\eqb{caseBigger}
\|u_k - u\opt\|^2_{M_k} \ge (1-2\phi_k) \| u\kp - u\opt\|^2_{M\kp}.
\eqe   
Note that $\lim_{k\to\infty}\phi_i=0,$ and so we may assume without loss of generality that $1-2\phi_k>0$  (this assumption is only violated for finitely many $k$).

  Now, consider the case that \eqref{upper:cond} does not hold. We have
   \begin{align}
   \|u\kp-&u\opt \|^2_{M_k}  \nonumber
      \ge \|u\kp-u\opt \|^2_{M\kp} \\
      &-  \frac{\phi_k }{\tau\kp}\|x\kp-x\opt\|^2 - \frac{\phi_k }{\sigma\kp}\|y\kp-y\opt\|^2  \\
       \ge&  \|u\kp-u\opt\|_{M\kp}^2 - \phi_k C_x. \label{smallerThanC}
    \end{align}
     Applying \eqref{smallerThanC} to Lemma \ref{lemma:ineq} yields
  \eqb{caseSmaller}\|u_k - u\opt\|^2_{M_k} \ge \| u\kp - u\opt\|^2_{M\kp}- \phi_k C_x. 
  \eqe
From \eqref{caseBigger} and \eqref{caseSmaller}, it follows by induction that 
 \begin{align}\label{combinedBound}
 \|u_0 - u\opt\|^2_{M_0} \ge&  \prod_{i\in I_C} (1-2\phi_i)\| u\kp \\
    &- u\opt\|^2_{M\kp}- \sum_{i } \phi_i C_x 
 \end{align}
  where $ I_C= \{i \, | \, \frac{1}{\tau\kp}\|x\kp-x\opt\|^2 + \frac{1}{\sigma\kp}\|y\kp-y\opt\|^2 \ge C_x \}.$  Note again that we have assumed without loss of generality that $k$ is large, and thus  $1-2\phi_k>0.$

We can rearrange \eqref{combinedBound} to obtain
  $$  \| u\kp - u\opt\|^2_{M\kp}
\le \frac{\|u_0 - u\opt\|^2_{M_0}+ C_x \sum_{i } \phi_i }{ \prod_i (1-2\phi_i)}  <\infty$$
which shows that  $\|u_k - u\opt\|^2_{M_k}$ remains bounded.  

Finally, note that since $\{\tau_k\},$ $\{\sigma_k\}, $ and $\|u_k-u\opt \|_{M_k}$ are bounded from above, it follows  from \eqref{bound:upper} that $ \frac{1}{\tau_k}\|x_k-x\opt\|^2$ is bounded from above.  But $ \frac{1}{\sigma_k}\|y_k-y\opt\|^2$ is also bounded from above, and so $\| u_k-u\opt\|_{H_k}$ is bounded as well.
  \end{IEEEproof}

Lemma \ref{lemma:upper} established upper bounds on the sequence of iterates.  To complete our convergence proof we also need lower bounds on $ \|u_k-u\opt\|^2_{M_k}.$  Note that in the case of indefinite $M_k,$ this quantity may be negative.  The following result show that $ \|u_k-u\opt\|^2_{M_k}$ does not approach negative infinity.

\begin{lemma} \label{lemma:lower}
Suppose  the step sizes for Algorithm \ref{alg:pdhg} satisfy conditions \A{}, \B{}, and \C{}.  Then 
$$\|u\kp-u\opt\|^2_{M_k}\ge C_L$$
for some lower bound $C_L.$
\end{lemma}
\begin{IEEEproof}  
  In the case that  \Cone{} holds, Lemma \ref{lemma:posdef} tells us that $M_k$ is positive definite.  In this case $ \|\cdot \|^2_{M_k}$ is a proper norm and
   $$\|u\kp-u\opt\|^2_{M_k}\ge 0.$$
   
     In the case that \Ctwo{} holds, we have that either $X$ or $Y$ is bounded.  Assume without loss of generality that $Y$ is bounded.  We then have $\|y\|<C_y$ for all $y\in Y$.  We can then obtain
     \begin{align} \|u\kp-u&\opt \|_{M_k}^2 \ge  \frac{1}{\tau_k}\|x\kp-x\opt\|^2\\
       & - 2C_y\|A\|_{op} \|x\kp-x\opt\|  + \frac{4C_y^2}{\sigma_k}\nonumber \\
     \ge&  \frac{1}{\min\{\tau_k\}}\|x\kp-x\opt\|^2 \\
     &- 2C_y\|A\|_{op} \|x\kp-x\opt\|  + \frac{4C_y^2}{\min\{\sigma_k\}}. \label{bound:lower} 
  \end{align}
  Note that \eqref{bound:lower} is quadratic in $\|x\kp-x\opt\|$, and so this quantity is bounded from below.
  
\end{IEEEproof}

We now present one final lemma, which bounds a sum involving the differences between adjacent iterates.

\begin{lemma} \label{lemma:sum}
Under conditions \A{}, \B{}, and \C{}, we have 
$$  \sum_{k=1}^{n}\|u_k - u\|^2_{M_k}-\|u_k - u\|^2_{M_{k-1}} <C_\phi C_U+ C_\phi C_H \| u-u\opt \|^2 $$
where $C_H$ is a constant such that $\| u-u\opt \|^2_{H_k} < C_H\| u-u\opt \|^2.$
\end{lemma}
\begin{IEEEproof}  
 We expand the summation on the left side of \eqref{converge:ineq2} using the definition of $M_k$ to obtain
\begin{align} \begin{split} \label{converge:ineq4}
 \sum_{k=1}^{n}&\|u_k - u\|^2_{M_k}-\|u_k - u\|^2_{M_{k-1}} \\
   &=  \sum_{k=1}^{n}(\frac{1}{\tau_k}-\frac{1}{\tau_{k-1}})\|x_k-x\|^2 +(\frac{1}{\sigma_k}-\frac{1}{\sigma_{k-1}})\|y_k-y\|^2\\
      &\le \sum_{k=1}^{n}(1-\delta_{k-1})   \left(\frac{1}{\tau_k}\|x_k-x\|^2 +\frac{1}{\sigma_k}\|y_k-y\|^2 \right)\\
    &=  \sum_{k=1}^{n}\phi_{k-1} \| u_k-u \|^2_{H_k}\\ 
     &\le  \sum_{k=1}^{n}\phi_{k-1} \left ( \| u_k-u\opt \|^2_{H_k} + \| u-u\opt \|^2_{H_k} \right)\\
    &   \le  \sum_{k=1}^{n}\phi_{k-1} \left ( C_U+ C_H\| u-u\opt \|^2 \right )\\
    & \le C_\phi C_U+ C_\phi C_H \| u-u\opt \|^2 <\infty, 
    \end{split}\end{align}
   where we have used the bound $\| u_k-u\opt \|^2_{H_k}<C_U$ from Lemma \ref{lemma:upper}.
    
\end{IEEEproof}

  \subsection{Convergence Theorems} \label{sec:thrms}
In this section, we prove convergence of Algorithm \ref{alg:pdhg} under assumptions \A{}, \B{}, and \C{}.
 The idea is to show that the norms of the residuals have a finite sum and thus converge to zero.  Because the ``natural'' norm of the problem, the $M_k$-norm, changes after each iteration, we must be careful to bound the differences between the norms used at each iteration.  For this purpose, condition \B{} will be useful because it guarantees that the various $M_k$-norms do not differ too much as $k$ gets large.

\begin{theorem} \label{backtrackConverge}
Suppose that the stepsizes in Algorithm \ref{alg:pdhg} satisfy conditions \A{} and \B{}, and either \Cone{} or \Ctwo{}.
Then the algorithm converges in the residuals, i.e.
    $$\lim_{k\to \infty} \|P_k\|^2 + \|D_k\|^2 = 0.$$
\end{theorem}
\begin{IEEEproof}

Rearranging Lemma \eqref{lemma:ineq} gives us
\begin{align} \label{converge:ineq2}
\|u_k - u\opt\|^2_{M_k} -   \| u\kp - u\opt\|^2_{M_k} \ge \| u\kp - u_k\|^2_{M_k}.
\end{align}
Summing \eqref{converge:ineq2} for $1\le k\le n$ gives us
\begin{align}\begin{split} \label{converge:ineq3}
\sum_{k=1}^{n} \| u\kp - &u_k\|^2_{M_k} \le \|u_1 - u\opt\|^2_{M_{0}}-\|u_{n+1} - u\opt\|^2_{M_{n}}\\
    & + \sum_{k=1}^{n}\|u_k - u\opt\|^2_{M_k}-\|u_k - u\opt\|^2_{M_{k-1}}.
\end{split}\end{align}

Applying Lemma \ref{lemma:sum} to \eqref{converge:ineq3}, and noting from Lemma \ref{lemma:lower} that $\|u_{n+1}-u\opt\|_{M_{n}}$ is bounded from below,  we see that
 $$\sum_{k=1}^{n} \| u\kp - u_k\|^2_{M_k}<\infty.$$
It follows that $\lim_{k\to \infty}\|u\kp-u_k\|_{M_k}=0.$ 
  If condition \Ctwo{} holds, then this clearly implies that 
  \eqb{hToZero}\lim_{k\to \infty}\|u\kp-u_k\|_{H_k}=0.\eqe
  If condition \Cone{} holds,  then we still have \eqref{hToZero} from Lemma \ref{lemma:posdef}.
  Since $\tau_k$ and $\sigma_k$ are bounded from above,  \eqref{hToZero} implies that
   \eqb{twoToZero}\lim_{k\to \infty}\|u\kp-u_k\|=0.\eqe
  
 Recall the definition of $R_k$ in equation \eqref{vecs}.   From \eqref{vi}, we know that
 $$R_k = M_k(u_k-u\kp) \in R(u\kp),$$ and so 
 \begin{align}\label{converge:lims}
\lim_{k\to \infty}\|R_k\| &=  \lim_{k\to \infty}\|M(u\kp-u_k)\| \\
 &\le \max_k\{\|M_k\|\} \lim_{k\to \infty}\|u\kp-u_k\|=0.
 \end{align}
\end{IEEEproof}

\begin{theorem} \label{convergeRate}
Suppose that the stepsizes in Algorithm \ref{alg:pdhg} satisfy conditions \A{}, \B{}, and \C{}.
Consider the sequence defined by
$$\tilde u_t = \frac{1}{t}\sum_{k=1}^t u_k.$$
  This sequence satisfies the convergence bound
  \begin{align*}
&\phi(u)-\phi( \tilde u_t ) +( u- \tilde u_t)^T  Q(\tilde u_t) \ge \\ 
& \frac{\| u-u_t \|^2_{M_{t}} - \| u-u_0 \|^2_{M_0}  - C_\phi C_U- C_\phi C_H \| u-u\opt \|^2}{2t}.
\end{align*}
  \end{theorem}
\begin{IEEEproof}

 We begin with the following identity (a special case of the polar identity for normed vector spaces):
 \begin{multline*}
 (u-u\kp)^TM_k(u_k-u\kp) = \half (  \| u-u\kp \|^2_{M_k}  \\ -  \| u-u_k \|^2_{M_k}  ) +\half \| u_k-u\kp \|^2_{M_k}.  
 \end{multline*}
We apply this to the VI formulation of the PDHG iteration \eqref{VI}
 to get
\begin{multline}\label{expanded}
 \phi(u)-\phi(u\kp)+(u-u\kp)^T  Q(u\kp)   \\ \ge  \half \left(  \| u-u\kp \|^2_{M_k} -  \| u-u_k \|^2_{M_k}  \right) \\ +\half \| u_k-u\kp \|^2_{M_k}.
\end{multline}
Note that 
\begin{multline}
(u-u\kp)^T Q(u-u\kp)=(x-x\kp)A^T(y-y\kp)\\ - (y-y\kp)A(x-x\kp)=0,
\end{multline} 
and so
\eqb{rateEqual}
(u-u\kp)^T  Q(u) = (u-u\kp)^T  Q(u\kp).
\eqe
  Also, both conditions \C1 and \C2 guarantee that $$\| u_k-u\kp \|^2_{M_k} \ge 0.$$  These observations reduce \eqref{expanded} to 
\begin{multline}\label{reduced}
 \phi(u)-\phi(u\kp)+(u-u\kp)^T  Q(u)  \\ \ge  \half \left(  \| u-u\kp \|^2_{M_k} -  \| u-u_k \|^2_{M_k}  \right) .
 \end{multline}
We now sum \eqref{reduced} for $k=0$ to $t-1,$ and invoke Lemma \ref{lemma:sum}.
\begin{align}\begin{split} \label{rateLower}
2\sum_{k=0}^{t-1}& \phi(u)-\phi(u\kp)+(u-u\kp)^T  Q(u)  \\ 
  \ge&  \| u-u_t \|^2_{M_{t}} - \| u-u_0 \|^2_{M_0}  \\
   &+ \sum_{k=1}^t   \left(  \| u-u_k \|^2_{M_{k-1}} -  \| u-u_k \|^2_{M_k}  \right) \\
\ge&  \| u-u_t \|^2_{M_{t}} - \| u-u_0 \|^2_{M_0} \\
 & - C_\phi C_U- C_\phi C_H \| u-u\opt \|^2.
\end{split}\end{align}
Because $\phi$ is convex, $$\sum_{k=0}^{t-1}\phi(u\kp)\le t\phi\left(\frac{1}{t}\sum_{k=1}^{t}  u_k \right)=t\phi(\tilde u_t).$$  The left side of \eqref{rateLower} therefore satisfies
\begin{multline}\label{rateConvex}
 2t \left( \phi(u)-\phi( \tilde u_t ) +( u- \tilde u_t)^T  Q(u) \right) \\ \ge 2\sum_{k=0}^{t-1} \phi(u)-\phi(u\kp)+(u-u\kp)^T  Q(u) .
\end{multline}
Combining  \eqref{rateLower} and \eqref{rateConvex} yields the tasty bound
\begin{align*}\begin{split} 
&\phi(u)-\phi( \tilde u_t ) +( u- \tilde u_t)^T  Q(u)  \\ 
&\ge   \frac{\| u-u_t \|^2_{M_{t}} - \| u-u_0 \|^2_{M_0}  - C_\phi C_U- C_\phi C_H \| u-u\opt \|^2}{2t}.
\end{split}\end{align*}
Applying \eqref{rateEqual} proves the theorem.
 \end{IEEEproof}

Note that  Lemma \ref{lemma:upper} guarantees $\{u_k\}$ remains bounded, and thus $\{\tilde u_t\}$ is bounded also.  Furthermore, because $\{\tau_k\}$ and $\{\sigma_k\}$ remain bounded, the spectral radii of the matrices $\{M_t\}$ remain bounded as well.   It follows that, provided $\|u-u_t\|\le1,$ there is a $C$ with
$$ \| u-u_0 \|^2_{M_0} -\| u-u_t \|^2_{M_{t}} \le C  \quad \forall\, t>0, u\in B_1(\tilde u_t).$$
We then have
 \begin{align*}
\phi(u)&-\phi( \tilde u_t ) +( u- \tilde u_t)^T  Q(\tilde u_t)  \\ 
&\ge \frac{-C - C_\phi C_U- C_\phi C_H \| u-u\opt \|^2}{2t}
\end{align*}
and the algorithm converges with rate $O(1/t)$ in an ergodic sense.

\section{Numerical Results}

 To demonstrate the performance of the new adaptive PDHG schemes, we apply them to the test problems described in Section \ref{section:problems}.  We run the algorithms with parameters $\alpha_0 = 0.5,$ $\Delta = 1.5,$ and $\eta = 0.95.$  The backtracking method is run with $\gamma = 0.75, \beta = 0.95,$ and we initialize stepsizes using formula \eqref{init}.   We terminate the algorithms when both the primal and dual residual norms (i.e. $|P_k|$ and $|D_k|$) are smaller than 0.05, unless otherwise specified.
 
 We consider four variants of PDHG. The method ``Adapt:Backtrack'' denotes Algorithm \ref{alg:adapt} with backtracking added. The initial stepsizes are chosen according to \eqref{init}.  The method ``Adapt: $\tau\sigma=L$'' refers to the adaptive method without backtracking with $\tau_0=\sigma_0=0.95\rho(A^TA)^{-\half}.$
 
  We also consider the non-adaptive PDHG with two different stepsize choices.  The method ``Const: $\tau,\sigma=\sqrt{L}$'' refers to the constant-stepsize method with both stepsize parameters equal to $\sqrt{L}=\rho(A^TA)^{-\half}. $  The method ``Const: $\tau$-final'' refers to the constant-stepsize method, where the stepsizes are chosen to be the final values of the stepsizes used by ``Adapt: $\tau\sigma=L$''.  This final method is meant to demonstrate the performance on PDHG with a stepsize that is customized to the problem at hand, but still non-adaptive.  
 
 \subsection{Experiments}
 The specifics of each test problem are described below:
 
 \subsubsection*{Rudin-Osher-Fatemi  Denoising}
  We apply the denoising model \eqref{rof} to the ``Cameraman'' test image.  The image is scaled to have pixels in the range $[0,255],$ and contaminated  with Gaussian noise of standard deviation 10.  The image is denoised with $\mu=0.25,$ $0.05,$ and $0.01.$  
     \begin{figure*}
     \centering
  \includegraphics[trim = .5cm 2cm 1cm .5cm, clip, width = 8.5cm]{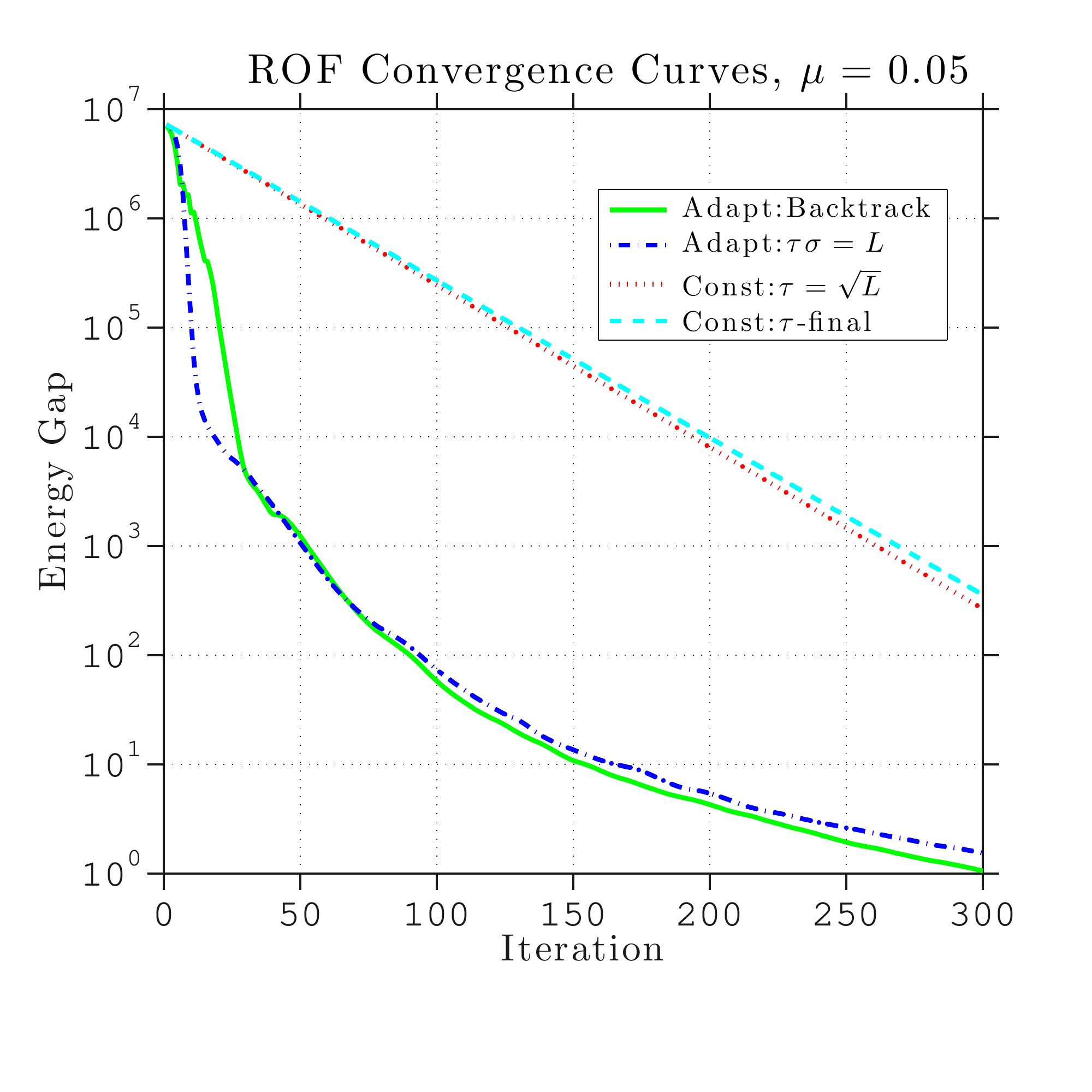}
   \includegraphics[trim = .5cm 2cm 1cm .5cm, clip, width = 8.5cm]{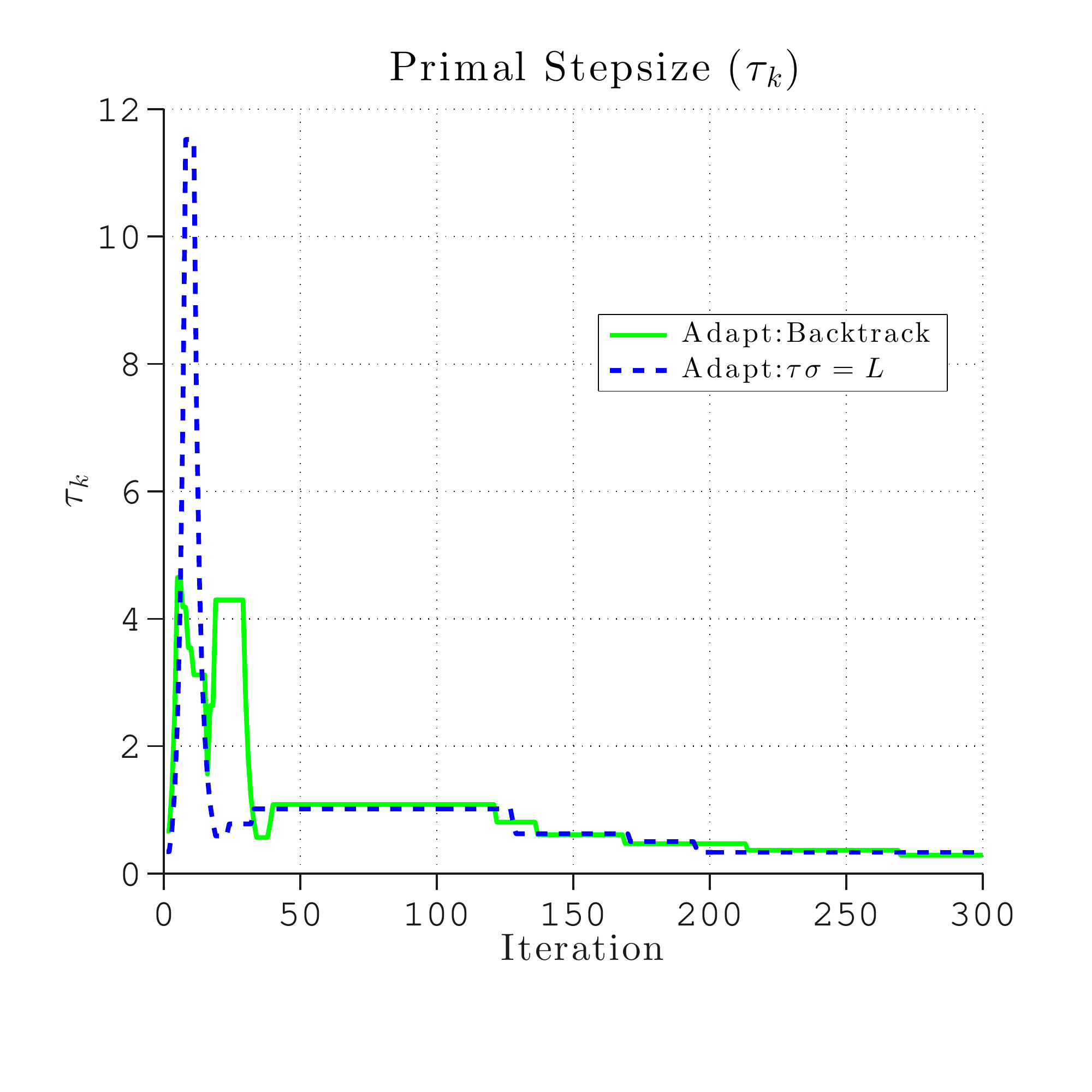}\\
         \caption{(left) Convergence curves for the Rudin-Osher-Fatemi denoising experiment with $\mu=0.05.$ The $y$-axis displays the difference between the value of the ROF objective function \eqref{rof} at the $k$th iterate and the optimal objective value.  (right) Stepsize sequences, $\{\tau_k\},$ for both adaptive schemes.}
    \label{fig:taus}
  \end{figure*}
   
   We display denoised images in Figure \ref{fig:denoise}. We show  results of numerical time trials in Table \ref{table}.  Note that the iteration count for denoising problems increases for small $\mu,$ which results in solutions with large piecewise -constant regions.  Note also the similar performance of Algorithm \ref{alg:adapt} with and without backtracking, indicating that there is no significant advantage to knowing the constant $L=\rho(A^TA)^{-1}.$
   
     We plot convergence curves and show the evolution of $\tau_k$ in Figure \ref{fig:taus}.  Note that $\tau_k$ is large for the first several iterates and then decays over time.  This behavior is typical for many TV-regularized problems.
   
\subsubsection*{TVL1 Denoising}
 We again denoise the ``Cameraman'' test image, this time using the model  \eqref{tvl1}, which tends to result in smoother results.  The image is denoised with $\mu=2,$ $1,$ and $0.5$.
 
\begin{figure}
  \centering
  \center{ \hspace{.5cm}\large{ROF}  \hspace{3.9cm} \large{TVL1}} \vspace{0.1cm} \\
  \includegraphics[trim = 1.5cm 1.5cm 1.5cm 1.5cm,clip, width = 4cm,height = 3.5cm]{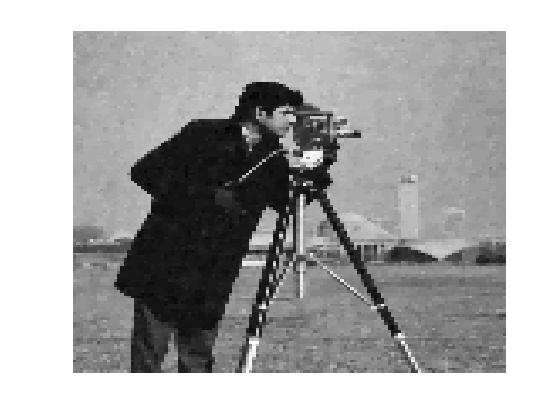}
   \includegraphics[trim = 1.5cm 1.5cm 1.5cm 1.5cm,clip, width = 4cm,height = 3.5cm]{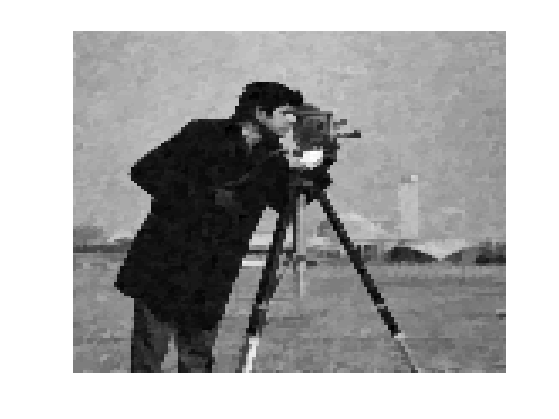}\\
    \vspace{-.1cm}$\mu=0.25$\hspace{3.9cm} $\mu=2$  \vspace{.1cm}\\
    \includegraphics[trim = 1.5cm 1.5cm 1.5cm 1.5cm,clip, width = 4cm,height = 3.5cm]{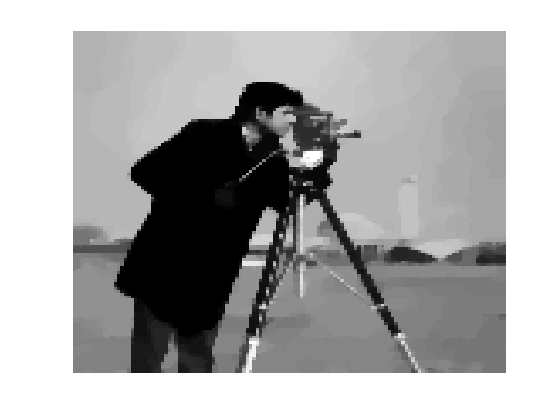}
     \includegraphics[trim = 1.5cm 1.5cm 1.5cm 1.5cm,clip, width = 4cm,height = 3.5cm]{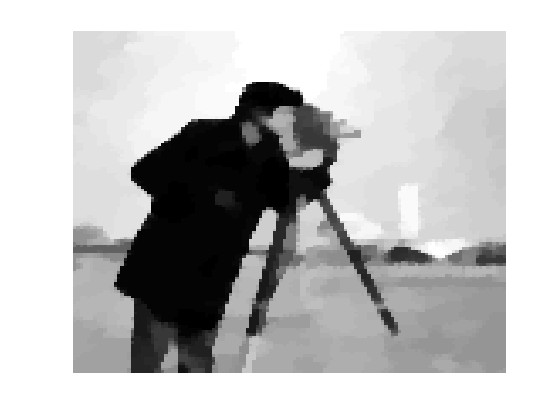}\\  
     \vspace{-.1cm} $\mu=0.05$\hspace{3.9cm} $\mu=1$ \vspace{.1cm}\\
     \includegraphics[trim = 1.5cm 1.5cm 1.5cm 1.5cm,clip, width = 4cm,height = 3.5cm]{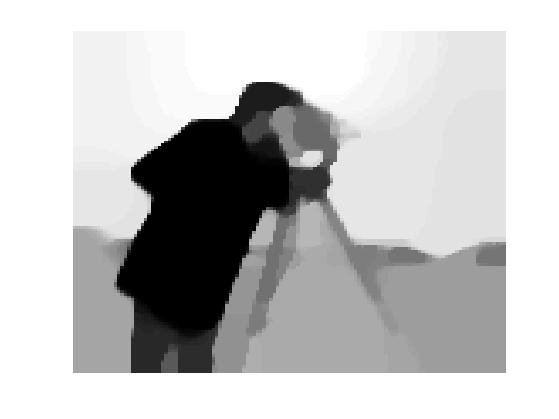}
      \includegraphics[trim = 1.5cm 1.5cm 1.5cm 1.5cm,clip, width = 4cm,height = 3.5cm]{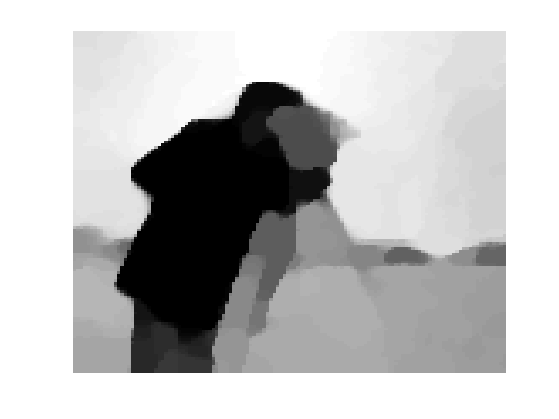}\\
      \vspace{-.1cm} $\mu=0.01$ \hspace{3.9cm} $\mu=0.5$ \vspace{.1cm}
        \caption{Results of denoising experiments with cameraman image.  (left column) ROF results with $\mu = 0.25,$ $0.05,$ and  $0.01$ from top to bottom.  (right column) TVL1 results with $\mu = 2,$ $1,$ and  $0.5$ from top to bottom.}
          \label{fig:denoise}
  \end{figure}

 We display denoised images in Figure \ref{fig:denoise}, and show time trials results in Table \ref{table}.  Much like in the ROF case, the iteration counts increase as denoising results get coarser (i.e. when $\mu$ gets small.)  There is no significant advantage to specifying the value of $L=\rho(A^TA)^T,$ because the backtracking algorithm was very effective for this problem, much like in the ROF case.

\subsubsection*{Convex Segmentation}
\begin{figure*}
 \centering
  \includegraphics[trim = 2.5cm 2.5cm 2.5cm 2.5cm,clip, width = 4.5cm]{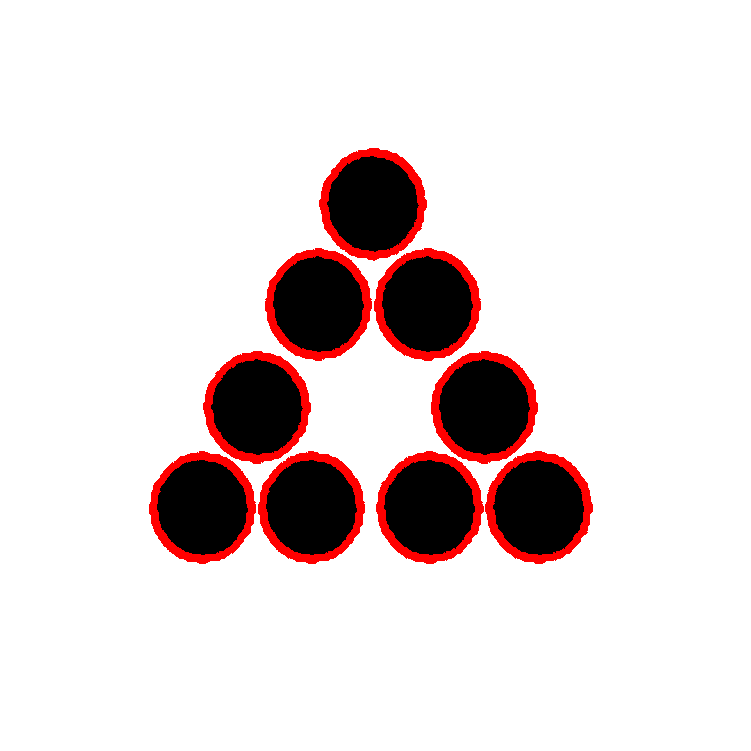}
   \includegraphics[trim = 2.5cm 2.5cm 2.5cm 2.5cm,clip, width = 4.5cm]{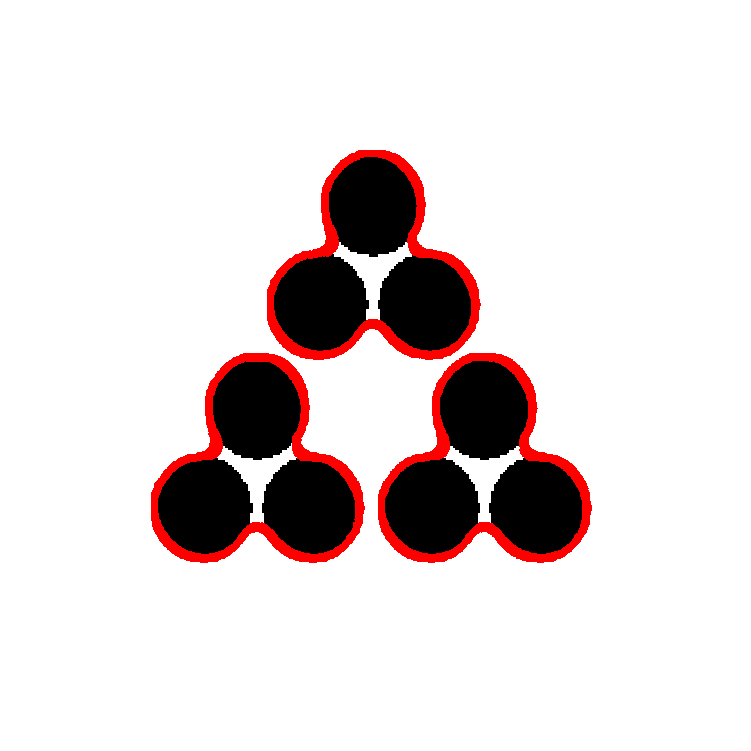}
    \includegraphics[trim = 2.5cm 2.5cm 2.5cm 2.5cm,clip, width = 4.5cm]{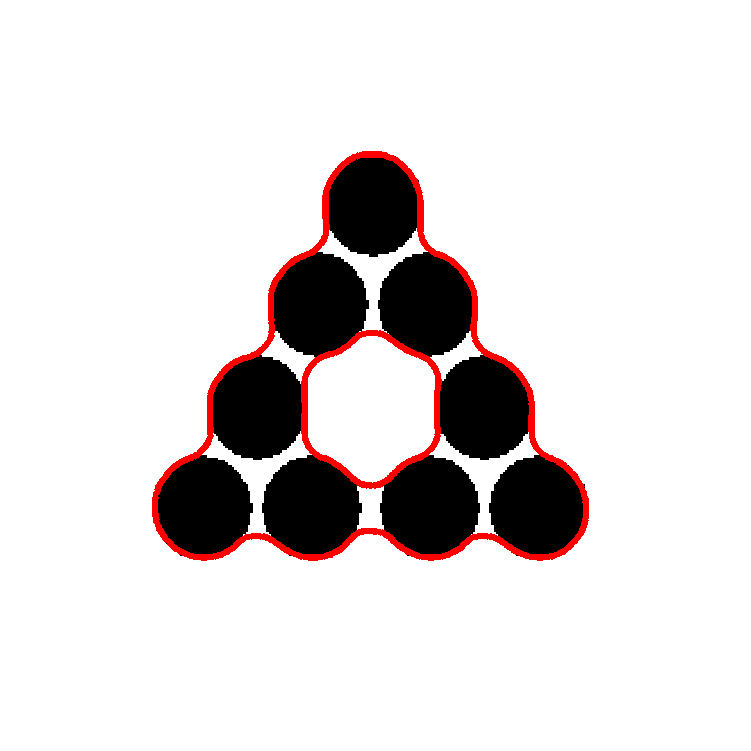}
           \caption{Segmentation of the ``circles'' test image at different scales.}
            \label{fig:seg}
  \end{figure*}

We apply the model \eqref{gcs} to a test image containing circular regions organized in a triangular pattern.  By choosing different weights for the data term $\mu$, we achieve  segmentations at different scales.  In this case, we can identify each circular region as its own entity, or we can groups regions together into groups of 3 or 9 circles.  
   Results of segmentations at different scales are presented in Figure \ref{fig:seg}.  

\subsubsection*{Compressed Sensing}
\begin{figure}
  \includegraphics[trim = 3.7cm 2cm 3cm 1cm, clip, width = 4.3cm,]{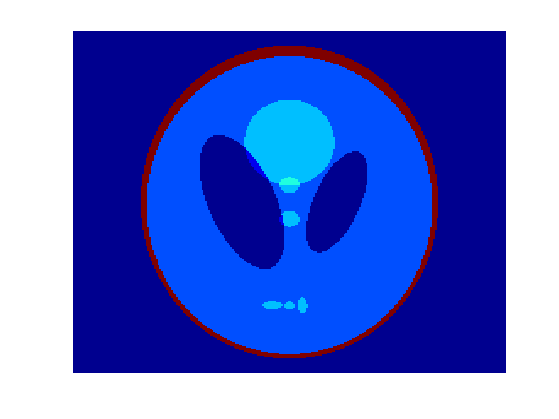}~
   \includegraphics[trim = 3.7cm 2cm 3cm 1cm, clip, width = 4.3cm]{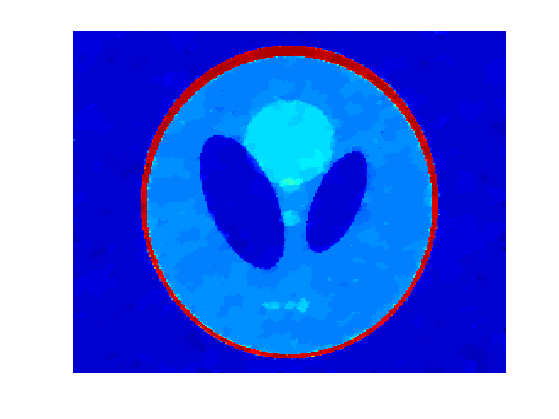}
      \caption{Compressed sensing reconstruction experiment.  (left) Original Shepp-Logan phantom. (right) Image reconstructed from a random 10\% sample of noisy Hadamard coefficients.  Images are depicted in false color to accentuate small features.}
        \label{fig:cs}
  \end{figure}

  We reconstruct a Shepp-Logan phantom from sub-sampled Hadamard measurements.  Data is generated by applying the Hadamard transform to a $256\times 256$ discretization of the Shepp-Logan phantom, and then sub-sampling $5\%,$ $10\%,$ and $20\%$ of the coefficients are random.  We scaled the image to have pixels in the range  $[0,255],$ the orthogonal scaling of the Hadamard transform was used, and we reconstruct all images with $\mu=1$.   The compressive reconstruction with $10\%$ sampling is shown in Figure \ref{fig:cs}.  See Table \ref{table} for iteration counts at the various sampling rates.
  
\subsubsection*{$\ell_\infty$ Minimization}   
To demonstrate the performance of the adaptive schemes with complex inputs, we solve a signal approximation problem using a complex-valued tight frame.  Problems of the form \eqref{linf} were generated by randomly sub-sampling 100 rows from a discrete Fourier matrix of dimension $512$.  The input signal $z$ was a vector on random Gaussian variables.  Problems were solved with approximation accuracy $\epsilon = 1,0.1,$ and $0.01,$ and results are reported in Table \ref{table}. 

  Because this problem is complex-valued, we use the definition of $b_k$ given by \eqref{backComp} to guarantee that the stepsizes are real-valued.
     
\subsubsection*{General Linear Programming}
We test our algorithm on the standard test problem ``sc50b'' that ships with most distributions of Matlab.  The problem is to recover 40 unknowns subject to 30 inequality and 20 equality constraints.  To accelerate the method, we apply a preconditioner to the problem inspired by the diagonal stepsize proposed for linear programming in \cite{PC11}.  This preconditioner works by replacing $A$ and $b$ in \eqref{lp} with 
 $$\hat A = \Gamma^\half A \Sigma^\half,\quad \hat b = \Gamma^\half b$$
 where $\Gamma$ and  $\Sigma$ are diagonal preconditioning matrices with $\Gamma_{ii} = \sum_{j}|A_{i,j}|$ and $\Sigma_{jj} = \sum_{i}|A_{i,j}|.$

Time trial results are shown in Table \ref{table}.  Note that PDHG is not as efficient as conventional linear programing methods for small-scale problems;  the Matlab ``linprog'' command took approximately 0.05 seconds to solve this problem using interior point methods.  Interestingly, the backtracking variant of Algorithm \ref{alg:adapt} out-performed the non-backtracking algorithm for this problem (see Table \ref{table}).
 
\begin{table*}[ht] \small
\caption{Time Trial Results.  Iteration counts are reported for each problem/algorithm, with total runtime (sec) in parenthesis.} 
\centering 
\begin{tabular}{c | c c c c} 
\hline\hline 
Problem & Adapt:Backtrack & Adapt: $\tau\sigma=L$ & Const: $\tau,\sigma=\sqrt{L}$ &  Const: $\tau$-final \\ [0.5ex] 
\hline 
ROF, $\mu=.25$ & 16 (0.0475)  & 16 (0.041) & 78 (0.184) & 48 (0.121) \\
ROF, $\mu=.05$  & 50 (0.122) & 51 (0.122) & 281 (0.669) & 97 (0.228) \\
ROF, $\mu=.01$ & 109 (0.262) & 122 (0.288) & 927 (2.17) & 152 (0.369) \\ 
TVL1, $\mu=2$ & 286 (0.957) & 285 (0.954) & 852 (2.84)  & 380 (1.27)\\ 
TVL1, $\mu=1$ & 523 (1.78) & 521 (1.73) & 1522 (5.066) & 669 (2.21) \\
TVL1, $\mu=.5$  & 846 (2.74) & 925 (3.07) & 3244 (10.8) & 1363 (4.55) \\
Segment, $\mu=0.5$ & 42 (0.420) & 41 (0.407) & 114 (1.13) & 53 (0.520) \\
Segment, $\mu=.15$ & 111 (1.13) & 107 (1.07) & 493 (5.03) & 131 (1.34) \\
Segment, $\mu=.08$ & 721 (7.30) & 1016 (10.3) & 706 (7.13) & 882 (8.93) \\
Compressive (20\%) & 163 (4.08)  & 168 (4.12) & 501 (12.54) & 246 (6.03) \\
Compressive (10\%) &  244 (5.63) & 274 (6.21) & 908 (20.6) & 437 (9.94) \\
Compressive \, (5\%) & 382 (9.54) & 438 (10.7) & 1505 (34.2) & 435 (9.95) \\
$\ell_\infty$ \, $(\epsilon=1)$ & 86 (0.0675) & 56 (0.0266) & 178 (0.0756) & 38 (0.0185) \\
$\ell_\infty$ \,$ (\epsilon=.1)$ & 79 (0.0360) &  72 (0.0309) & 195 (0.0812) & 78 (0.0336) \\
$\ell_\infty$ \, $(\epsilon=.01) $& 89 (0.0407) & 90 (0.0392) & 200 (0.0833) & 82 (0.0364) \\
LP & 18255 (3.98) & 20926 (4.67) & 24346 (5.42) & 29357 (6.56)\\
 [1ex] 
\hline 
\end{tabular}
\label{table} 
\end{table*}

\subsection{Discussion}  Several interesting observations can be made from the results in Table \ref{table}.  First, both the backtracking (``Adapt:Backtrack'') and non-backtracking (``Adapt: $\tau\sigma=L$'') methods have similar performance on average for the imaging problems, with neither algorithm showing consistently better performance than the other.  This shows that the backtracking scheme  is highly effective and that there is no significant advantage to requiring the user to supply the constant $\rho(A^TA)$  other than for simplicity of implementation.   

Interestingly, when the backtracking method did show better performance, it was often due to the ability to use large stepsizes that violate the stepsize restriction \Cone{}.  When solving the ROF denoising problem with $\mu=0.01$, for example, the backtracking scheme terminated with $\tau_k\sigma_k = 0.14,$ while the conventional stepsize restriction is $1/\rho(A^TA) = 1/8=0.125.$ 


Finally, the method  ``Const: $\tau$-final'' (a non-adaptive method using the ``optimized'' stepsizes obtained from the last iterations of the adaptive scheme) did not always outperform the non-adaptive scheme using the non-optimized stepsizes.  This occurs because the true ``best'' stepsize choice depends on the active set of the problem in addition to the structure of the remaining error and thus evolves over time.  This situation is depicted in Figure \ref{fig:taus}, which shows the evolution of $\tau_k$ over time for the adaptive methods.  Note that for this problem, both methods favor a large value for $\tau_k$ during the first few iterations, and this value decays over time.  Behavior like this is typical for the imaging problems discusses above.  This illustrates the importance of an adaptive stepsize --- by optimizing the stepsizes to the current active set and error levels, adaptive methods can achieve superior performance.

\section{Conclusion}
  We have introduced new adaptive variants of the powerful Primal-Dual Hybrid Gradient (PDHG) schemes.     
   We proved rigorous convergence results for adaptive methods and identified the necessary conditions for convergence.  We also presented practical implementations of adaptive schemes that formally satisfy the convergence requirements.  Numerical experiments show that adaptive PDHG methods have advantages over non-adaptive implementations in terms of both efficiency and simplicity for the user.  
    
    The new backtracking variant of PDHG is advantageous because it does not require the user to supply any information about the problem instance being solved such as the spectral radius of the matrix $A$.  This is particularly useful when creating ``black-box'' solvers for large classes of problems.  For example, a generic linear programming solver should not require the user to supply information about the eigenvalues of the constraint matrix.  Furthermore, this flexibility seems to come at no added cost;  the fully adaptive backtracking scheme performs at least as well as methods that do require spectral information.

\section{Acknowledgements}
 This work was made possible by the generous support of the CIA Postdoctoral Fellowship (\#2012-12062800003), and the NSF Post-doctoral fellowship program.

\bibliography{/Users/Tom/Documents/latexDocs/bib/tom_bibdesk}

\end{document}